\documentclass[11pt]{amsart}
\usepackage{amsxtra}
\usepackage{amssymb,yhmath,stmaryrd,multicol}
\theoremstyle{plain}

\newcommand{\cleqn}{\setcounter{equation}{0}}
\newcommand{\clth}{\setcounter{theorem}{0}}
\newcommand {\sectionnew}[1]{\section{#1}\cleqn\clth}

\newtheorem{theorem}{Theorem}[section]
\newtheorem{lemma}[theorem]{Lemma}
\newtheorem{definition-theorem}[theorem]{Definition-Theorem}
\newtheorem{proposition}[theorem]{Proposition}
\newtheorem{corollary}[theorem]{Corollary}
\newtheorem{definition}[theorem]{Definition}
\newtheorem{example}[theorem]{Example}
\newtheorem{remark}[theorem]{Remark}
\newtheorem{conjecture}[theorem]{Conjecture}
\newtheorem{notation}[theorem]{Notation}
\newcommand \bth[1] { \begin{theorem}\label{t#1} }
\newcommand \ble[1] { \begin{lemma}\label{l#1} }

\newcommand \bpr[1] { \begin{proposition}\label{p#1} }
\newcommand \bco[1] { \begin{corollary}\label{c#1} }
\newcommand \bde[1] { \begin{definition}\label{d#1}\rm }
\newcommand \bex[1] { \begin{example}\label{e#1}\rm }
\newcommand \bre[1] { \begin{remark}\label{r#1}\rm }
\newcommand \bcj[1] { \begin{conjecture}\label{j#1}\rm }

\newcommand \bnota[1] { \begin{notation}\label{n#1}\rm }
\renewcommand {\eth} { \end{theorem} }
\newcommand {\ele} { \end{lemma} }

\newcommand {\epr} { \end{proposition} }
\newcommand {\eco} { \end{corollary} }
\newcommand {\ede} { \end{definition} }
\newcommand {\eex} { \end{example} }
\newcommand {\ere} { \end{remark} }
\newcommand {\ecj} { \end{conjecture} }

\newcommand {\enota} { \end{notation} }
\newcommand \thref[1]{Theorem \ref{t#1}}
\newcommand \leref[1]{Lemma \ref{l#1}}
\newcommand \prref[1]{Proposition \ref{p#1}}

\newcommand \deref[1]{Definition \ref{d#1}}

\newcommand \reref[1]{Remark \ref{r#1}}
\newcommand \lb[1]{\label{#1}}

\def \Rset {{\mathbb R}}         

\def \KK {k}
\def \Zset {{\mathbb Z}}
\def \Nset {{\mathbb N}}

\def \Ksp {K^{\mathrm{split}}_0}
\def \AA  {{\mathcal{A}}}           

\def \FF {{\mathscr{F}}}

\def \PP {{\mathcal{P}}}
\def \MM {{\mathscr{M}}}

\def \cO {{\mathcal{O}}}

\def \De {\Delta}   

\def \al {\alpha}

\def \Sig {\Sigma}

\def \ep {\epsilon}

\def \mt  {\mapsto}
\def \lra {\longrightarrow}

\def \hra {\hookrightarrow}

\def \sup {\supset}

\def \rcor {\rangle}
\def \lcor {\langle}
\def \ol {\overline}

\def \id { {\mathrm{id}} }

\def \const { {\mathrm{const}} }

\def \add {{\mathrm{add}}}
\def \ind {{\mathrm{ind}}}
\def \Indec {{\mathrm{Indec}}}
\def \g  {\mathfrak{g}}   

\def \sl {\mathfrak{sl}} 
\def \hh  {\mathfrak{h}}

\def \b  {\mathfrak{b}}

\def \sl {\mathfrak{sl}}

\DeclareMathOperator \Span { {\mathrm{Span}} }
\DeclareMathOperator \Aut { {\mathrm{Aut}} }

\DeclareMathOperator \diag { {\mathrm{diag}} }

\DeclareMathOperator \Hom { {\mathrm{Hom}} }

\DeclareMathOperator \Ext { {\mathrm{Ext}} }

\renewcommand \Im { {\mathrm{Im}} }

\renewcommand \max { {\mathrm{max}} }

\usepackage{mathdots}

\usepackage{color}

\usepackage{pdfsync}

\usepackage{enumitem}

\usepackage{wasysym}

\usepackage{amssymb}

\usepackage{mathrsfs}

\DeclareMathAlphabet{\mathpzc}{OT1}{pzc}{m}{it}

\usepackage[all]{xy}

\usepackage{tikz}
\usetikzlibrary{arrows,decorations.pathmorphing,decorations.pathreplacing,positioning,shapes.geometric,shapes.misc,decorations.markings,decorations.fractals,calc,patterns}

\usepackage{float}

\usepackage[bottom]{footmisc}

\entrymodifiers={+!!<0pt,\fontdimen22\textfont2>}

\def\cC{\mathscr{C}}

\def\cO{\mathscr{O}}

\def\cT{\mathscr{T}}
\def\cU{\mathscr{U}}

\def\BN{\mathbb{N}}

\def\BR{\mathbb{R}}

\def\add{\operatorname{add}}

\def\adots{\mathinner{\mkern1mu\raise1.0pt\vbox{\kern7.0pt\hbox{.}}\mkern2mu\raise4.0pt\hbox{.}\mkern2mu\raise7.0pt\hbox{.}\mkern1mu}}

\def\ast{{\textstyle *}}
\def\astsmall{{\scriptstyle *}}
\def\Aut{\operatorname{Aut}}

\def\supp{\operatorname{supp}}

\def\dddots{\mathinner{\mkern1mu\raise10.0pt\vbox{\kern7.0pt\hbox{.}}\mkern2mu\raise5.3pt\hbox{.}\mkern2mu\raise1.0pt\hbox{.}\mkern1mu}}
\def\dddotssmall{\mathinner{\mkern1mu\raise7.0pt\vbox{\kern7.0pt\hbox{.}}\mkern-1mu\raise4pt\hbox{.}\mkern-1mu\raise1.0pt\hbox{.}\mkern1mu}}

\def\dim{\operatorname{dim}}

\def\ext{\operatorname{ext}}

\def\Ext{\operatorname{Ext}}

\def\intt{\operatorname{int}}

\def\Hom{\operatorname{Hom}}
\def\id{\operatorname{id}}
\def\Image{\operatorname{Im}}

\def\inf{\operatorname{inf}}

\def\mod{\mathsf{mod}}

\def\opp{\operatorname{op}}

\def\SL2{\operatorname{SL}_2}

\def\sup{\operatorname{sup}}

\begin{document}
\title[Categorical C-vectors and $\sl_\infty$]
{C-vectors of 2-Calabi--Yau categories and \\
Borel subalgebras of $\sl_\infty$}
\author[Peter J{\o}rgensen]{Peter J{\o}rgensen}
\address{School of Mathematics and Statistics,
Newcastle University, Newcastle upon Tyne NE1 7RU, United Kingdom}
\email{peter.jorgensen@ncl.ac.uk}
\urladdr{http://www.staff.ncl.ac.uk/peter.jorgensen}

\author[Milen Yakimov]{Milen Yakimov}
\address{
Department of Mathematics \\
Louisiana State University \\
Baton Rouge, LA 70803 \\
U.S.A.
}
\email{yakimov@math.lsu.edu}




\begin{abstract} 
We develop a general framework for $c$-vectors of 2-Calabi--Yau categories, which deals with cluster tilting subcategories that are not reachable from each other and contain infinitely many indecomposable objects.  It does not rely on iterative sequences of mutations. 

We prove a categorical (infinite-rank) generalization of the Nakanishi--Zelevinsky duality for $c$-vectors and establish 
two formulae for the effective computation of $c$-vectors -- one in terms of indices and the other in terms of dimension vectors 
for cluster tilted algebras.

In this framework, we construct a correspondence between the $c$-vectors of the cluster categories $\cC(A_\infty)$ of type $A_\infty$ 
due to Igusa--Todorov and the roots of the Borel subalgebras of $\sl_\infty$. Contrary to the finite dimensional case, 
the Borel subalgebras of $\sl_\infty$ are not conjugate to each other. On the categorical side, the cluster tilting subcategories of $\cC(A_\infty)$
exhibit different homological properties. The correspondence builds a bridge between the two classes of objects.
\end{abstract}
\maketitle
\sectionnew{Introduction}
\lb{intro}

\subsection{$c$-vectors of 2-Calabi--Yau categories}
\label{1.1}
One of the cornerstones of the theory of cluster algebras is their link to Lie algebras:
\bth{FZthm} {\em{(Fomin--Zelevinsky \cite{FZ0,FZ})}} Finite type cluster algebras are in bijection with the Cartan matrices of finite 
dimensional simple Lie algebras.
Let $A$ be such a cluster algebra and $B$ be the corresponding Cartan matrix. For each seed $S$ of  $A$ whose symmetrized mutation matrix 
equals $B$,
\[
C^+_S(A) = D_S(A) = \Delta^+(B)
\]
where $C^+_S(A)$ is the set of positive $c$-vectors of $A$ with respect to the seed $S$, 
$D_S(A)$ is the set of non-initial $d$-vectors with respect to $S$ and $\De^+(B)$ are the 
positive roots of the roots system with Cartan matrix $B$.
\eth
There has been a great deal of research on the extension of the correspondence between $c$-vectors of cluster algebras 
and roots of Kac--Moody algebras.  Nakanishi--Stella \cite{NS} did this for arbitrary seeds of finite type cluster algebras, and N\'ajera Ch\'avez
\cite{Ch1,Ch2} extended the correspondence to special seeds of acyclic cluster algebras. Various categorical interpretations of 
$g$-, $c$- and $d$-vectors were found in \cite{BMR,DK,DWZ,Tr} and many other papers.

Dehy--Keller \cite{DK} introduced a general notion of $g$-vectors of 2-Calabi--Yau (2-CY) categories $\cC$ in terms of indices 
with respect to cluster tilting subcategories. The radically new feature of 
the approach of \cite{DK} was that it allowed the simultaneous treatment of cluster tilting subcategories not reachable from 
each other, and did not require iterative mutations. In the special case of 2-CY categories $\cC$ with cluster structures 
in the sense of \cite{BIRS} and rigid objects that are reachable from the original cluster tilting subcategory, 
this definition recovers the combinatorial one of Fomin and Zelevinsky.
In this generality, further properties of indices/$g$-vectors were 
obtained in \cite{JP,P,P2}. 

In this paper we develop a general framework for $c$-vectors of 2-CY categories $\cC$ with respect to arbitrary cluster tilting 
subcategories $\cT$, based on the Dehy--Keller treatment of $g$-vectors. This approach deals with cluster tilting subcategories which are in 
general unreachable from each other, and does not perform (finite or infinite) sequences of mutations. We prove a general categorical duality of 
$c$-vectors which generalizes the combinatorial duality of Nakanishi and Zelevinsky \cite{NZ} to cluster tilting subcategories with infinitely many indecomposable objects, without assuming that such subcategories are reachable from each other by mutation. We also establish two formulae for the effective computation 
of $c$-vectors.
We propose a general program for decomposing sets of $c$-vectors and identifying the pieces with root systems of Kac-Moody algebras:
\begin{enumerate}
\item[(I)] The set $C^+_{\cT} (\cC)$ is equipped with a natural poset structure, and as a consequence, 
it equals the union of its maximal ideals. (Here and below by a maximal ideal of a poset, we mean an 
ideal that is maximal with respect to inclusion, but not necessarily proper, see \deref{R}.)
\item[(II)] Prove that each maximal ideal of $C^+_{\cT} (\cC)$ is isomorphic to the set of (Schur) roots of the Borel subalgebra 
of a Levi factor of the Kac--Moody algebra $\g(\ol{Q}_\cT)$ for the non-oriented graph underlying the quiver of the category $\cT$. 
\end{enumerate}
We realize the program for the cluster categories of type $A_\infty$ introduced by Igusa and Todorov, with respect to all their cluster tilting subcategories.
In part (II) we obtain a relation to the root systems  of the Borel subalgebras of $\sl_\infty$. The structure of the Borel subalgebras of $\sl_\infty$ 
and the related representations have been the subject of intensive research initiated by Dimitrov--Penkov \cite{DP1,DP2}.

Unlike the finite rank case, the Borel subalgebras of infinite rank Kac--Moody algebras $\g$ are not conjugate to 
each other. Their root systems behave in a drastically different way from the standard set of positive roots of $\g$, see \cite{DP1,DP2}.

The new features of our program to relate $c$-vectors of cluster categories and root systems, compared to the original theorem
of Fomin and Zelevinsky and its previous generalizations, are as follows:
\begin{enumerate}
\item[(a)] We consider $c$-vectors of unreachable cluster tilting subcategories of 2-CY categories, because the 
collections of $c$-vectors will otherwise miss big parts of a root system.
\item[(b)] We consider all Borel subalgebras of Kac--Moody algebras of infinite rank, because the (positive) 
root systems obtained in part (II) of the program are rarely generated by simple roots, while the standard set of positive roots of a Kac--Moody 
algebra always has this property.
\item[(c)] We decompose the collection of $c$-vectors into a union of maximal subsets, each of which comes from a single Levi factor of the Kac--Moody algebra $\g(\ol{Q}_\cT)$, instead of finding combinatorial patterns 
determining which roots of $\g(\ol{Q}_\cT)$ come from $c$-vectors as was previously done in concrete acyclic situations.
\end{enumerate} 
In one direction, the realization of the program for a 2-CY category $\cC$ classifies the cluster tilting subcategories of $\cC$ in terms of root systems 
because by \cite{DK} non-isomorphic indecomposable rigid objects of $\cC$ have different indices. In the other direction, the program produces
explicit categorifications of root systems of Borel subalgebras of Kac--Moody algebras. A remarkable feature of the correspondence is  that
\begin{enumerate}
\item[(d)] If we fix the category $\cC$, but vary the cluster tilting subcategory $\cT$ with respect to which indices are computed, then 
we obtain Borel subalgebras of Kac--Moody algebras whose categories $\cO$ are quite different from each other. At the same time, 
their roots are the $c$-vectors of exactly the same family of cluster tilting subcategories of $\cC$.
\end{enumerate}
This may be an indication of the existence of a translation principle between the characters of the modules of those categories,
because their Verma modules have characters equal to those of the symmetric algebras of the nilradicals of the Borel subalgebras.

Our program can be viewed as a construction of additive categorifications of the root systems of the Borel subalgebras of infinite rank Kac--Moody algebras 
(such as $\sl_\infty$) in terms of cluster tilting subcategories of 2-CY categories. We expect that this categorical setting will be helpful in 
addressing representation theoretic questions for the Borel subalgebras.
\subsection{Results on $c$-vectors of 2-CY categories}
\label{1.2}
Let $\KK$ be an algebraically closed field, $\cC$ a $\KK$-linear Hom-finite Krull--Schmidt triangulated category 
which is 2-CY with suspension functor $\Sigma$. Assume that $\cC$ has a cluster tilting subcategory $\cT$. Denote by 
$\Ksp(\cC)$ and $\Ksp(\cT)$ the split Grothendieck groups of $\cC$ and $\cT$. For an object $s \in \cC$, one considers its
index (or $g$-vector) with respect to $\cT$ 
\[
\ind_{\cT}(s) \in \Ksp(\cT),
\] 
see \cite{DK} and and \S \ref{2.1} for details. This induces a homomorphism 
\[
\ind_{\cT} \colon \Ksp(\cC) \to \Ksp(\cT) \quad \mbox{such that} \quad \ind_{\cT}|_{\Ksp(\cT)} = \id_{\Ksp(\cT)}.
\]
Dehy--Keller \cite{DK} proved that for every other cluster tilting subcategory $\cU$ of $\cC$,
\[
\{ \ind_{\cT}(u) \mid u \in \Indec(\cU) \} \; \; \mbox{is a basis of} \; \; \Ksp(\cT),
\]
where $\Indec(\cU)$ denotes the set of indecomposable objects of $\cU$. Define the $c$-vector of the pair $(u,\cU)$ 
where $u \in \Indec(\cU)$ to be the unique
\[
c_{\cT}(u,\cU) \in \Ksp(\cT)^* \quad \mbox{such that} \quad
\lcor c_{\cT}(u,\cU), \ind_{\cT}(v) \rcor = \delta_{uv}, \; \forall v \in \Indec(\cU). 
\]
This generalizes the Fomin--Zelevinsky definition \cite{FZ}: for every 2-CY category $\cC$ with cluster structure \cite{BIRS}
and a pair of cluster tilting subcategories $\cT$ and $\cU$ which are reachable from each other, this definition 
recovers the cluster algebra $c$-vector of the cluster variable corresponding to $u$ in the seed corresponding to $\cU$, 
see \S \ref{2.2}. From \cite{DK} we derive that all $c$-vectors $c_{\cT}(u,\cU)$ are sign coherent.  Note that $c_{ \cT }( u,\cU )$ is non-zero, see Remark \ref{rnon-zero}.  

Our first result is a categorical duality for the index maps with respect to cluster tilting subcategories allowed to be unreachable from each other.
It generalizes the Nakanishi--Zelevinsky duality for $c$-vectors of cluster algebras \cite[Theorem 1.2, Eq. (1.12)]{NZ} 
to arbitrary 2-CY categories.  For cluster tilting subcategories with finitely many indecomposables, the result was proved by Demonet--Iyama--Jasso \cite[thm.\ 6.19(a)]{DIJ} 
using $\tau$-tilting theory.  We will give a direct proof in the general case.

\bth{1} For every 2-CY category $\cC$ (which is $\KK$-linear, Hom-finite, Krull--Schmidt) and every pair of cluster tilting subcategories
$\cT$ and $\cU$, there are inverse isomorphisms
\[
  \xymatrix {
  \Ksp( \cT )
    \ar[rrr]<1ex>^-{ \overline{\ind}_{ \cU }\big|_{ \Ksp( \cT ) } } & & &
    \Ksp( \cU ),
    \ar[lll]<1ex>^-{ \ind_{ \cT } \big|_{ \Ksp( \cU ) } }
            }
\]
where $\overline{ \ind }_{\cU}$ denotes the index map in the opposite category of $\cC$.
\eth

Our next result effectively computes $c$-vectors in two different ways. For $t \in \Indec(\cT)$, let 
\[
[t]^* \in \Ksp(\cT)^* \; \; \mbox{be such that} \; \; \lcor [t]^* , [s] \rcor = \delta_{ts}, \forall s \in \Indec(\cT)
\]
and let $[\Sigma \cT]$ be the ideal of morphisms in $\cC$ which factor through an object of $\Sigma \cT$.  
\bth{2} Assume the setting of \thref{1}. Let $u \in \Indec(\cU)$. 
\begin{enumerate}
\item
The $c$-vector of the 
pair $(u,\cU)$ with respect to the cluster tilting subcatgory $\cT$ is given by
\[
c_{\cT}(u, \cU) = [u]^* \circ \overline{ \ind }_{\cU} |_{\Ksp(\cT)}.
\]
\item Assume in addition that the quiver of each cluster tilting subcategory is without loops and $2$-cycles.  Denote by $u^*$ the mutation of $u$ with respect to $\cU$ and let $u^{ \ast } \rightarrow e \rightarrow u \xrightarrow{ \delta } \Sigma u^{ \ast }$ be one of the two exchange triangles. The $c$-vector $c_{\cT}(u, \cU)$ is positive if and only if $\delta \not\in [\Sigma \cT]$, and in this case
\[
c_{\cT}(u, \cU)
= \dim_{ \cT } \Image \Big( \Hom(-,u)\big|_{\cT} \rightarrow \Hom (-,\Sigma u^*)\big|_{\cT} \big). 
\]
\end{enumerate}
\eth
Both parts of the theorem provide effective ways to evaluate $c$-vectors, in that no infinite matrices are inverted as in the original
definition of $c$-vectors. The second part generalizes a theorem of N\'ajera Ch\'avez \cite{Ch2}, which deals with the case of 
Amiot's cluster categories and a pair of cluster tilting subcategories which are reachable from each other. 

We obtain further properties of the set of $c$-vectors $C_{\cT}(\cC)$ of a 2-CY category $\cC$ with respect to a given cluster tilting subcategory $\cT$. We show
that $C_{\cT}(\cC)= - C_{\cT}(\cC)$ if $\cC$ has a cluster structure, and that the sets of $c$-vectors of $\cC$ and the dual category coincide.
We also show the needed decomposition of the set of positive $c$-vectors $C_{\cT}^+(\cC)$ for part (I) of the program in \S \ref{1.1}.
\subsection{Borel subalgebras of $\sl_\infty$ and $c$-vectors of cluster categories of type $A_\infty$}
\label{1.3}
Set $\Nset =\{0,1,\ldots\}$. The Lie algebra $\sl_\infty:=\sl_\infty(\KK)$ is the Lie algebra of $\Zset \times \Zset$ (or, equivalently, $\Nset \times \Nset$)
traceless matrices with only finitely many nonzero entries in a field $\KK$. Its standard Cartan subalgebra $\hh$ consists of diagonal matrices. The Borel subalgebras 
of $\sl_\infty$, containing $\hh$, are classified by the countable totally ordered sets \cite{DP1}.
The root system of the Borel subalgebra corresponding to such a set $Y$ is the subset
\[
\De^+_{\sl_\infty}(Y) := \{ \ep_{y'} - \ep_{y''} \mid y' > y'' \in Y \} 
\] 
of the root system of $\sl_\infty$ identified with $\De_{\sl_\infty} \cong \{ \ep_{y'} - \ep_{y''} \mid y' \neq y'' \in Y \}.$

The cluster categories $\cC(Z)$ of type $A_{ \infty }$ are 2-CY categories defined starting from an infinite subset $Z$ of the unit circle $S^1$, 
such that no point of $Z$ is a limit point of $Z$ and every limit point of $Z$ is both a left and a right limit point. The 
indecomposable objects of $\cC(Z)$ are classified by the diagonals $\{e,f\}$ of $Z$ (open segments joining non-consecutive vertices of $Z$). 
It was proved in \cite{GHJ,SvR} 
that the cluster tilting subcategories of $\cC(Z)$ are classified by the triangulations $T$ of $Z$ such that each accumulation 
point of $Z$ is surrounded by a fountain or a leapfrog.  Denote the corresponding category by $\cT(T)$.

In Sect. \ref{ind} we obtain explicit formulae for the $g$- and $c$-vectors of the categories $\cC(Z)$ with respect to any cluster tilting subcategory. 
Based on these we prove the following result linking $\cC(Z)$ to the Borel subalgebras of $\sl_\infty$. 

We will say that a totally ordered set $Y$ is {\em sequential} 
if every element of $Y$ has an immediate predecessor and an immediate successor, except, respectively,  
its least and greatest elements, if they exist.
Recall that a totally ordered set $Y$ is 
equipped with the order topology generated by the sets 
\[
(- \infty, y ) = \{ x \in Y \mid x< y\} \quad \mbox{and} \quad (y, + \infty) = \{ x \in Y \mid x > y \} \quad \mbox{for} \; \; y \in Y.
\] 
A totally ordered set $Y$ is sequential if and only if it is discrete in the order topology.

For a totally ordered set $Y$, 
denote the extended totally ordered set 
\begin{equation}
\label{extended}
Y_{\ext} := 
\begin{cases}
Y \sqcup \{ - \infty \}, & \; \; \mbox{if $Y$ has a least element} \\
Y,               & \; \; \mbox{otherwise}
\end{cases}
\end{equation}
with the relation $- \infty <y$ for all $y \in Y$.

\bth{3} Let $\cT(T)$ be an arbitrary cluster tilting subcategory of a cluster category $\cC(Z)$ of type $A_\infty$.
\begin{enumerate}
\item The maximal ideals of $C^+_{\cT(T)}(\cC(Z))$ are classified 
by the diagonals $\{e,f\}$ of $\ol{Z}$, with each $e$ and $f$ being an ear of $T$ or a limit point of a leapfrog of $T$.  
The corresponding maximal ideal $C^+_{\cT(T)}(\cC(Z))^\max_{e,f}$ consists of those $c$-vectors which are contained in
\[
\prod \{ \Zset \cdot [t]^* \mid t \in T, t \cap \{e,f\} \neq \varnothing \} \subset \Ksp(\cT(T))^*.
\]
\item If the dual quiver $Q_T$ of the triangulation $T$ has no cycles, then $C^+_{\cT(T)}(\cC(Z))$ has a unique maximal ideal 
(equal to itself).
\item The maximal ideal $C^+_{\cT(T)}(\cC(Z))^\max_{e,f}$ is (additively) isomorphic to 
\begin{itemize}
\item[(a)] the set of roots $\De^+_{\sl_\infty}(Y_{\ext})$ of the Borel subalgebra of $\sl_\infty$ corresponding to $Y_{\ext}$, where 
$Y:= \{e,f\} \cap T$ is ordered from $e$ to $f$ if $|Y|= \infty$, and
\item[(b)] the set $\Delta^+_{\sl_n}$ of positive roots of $\sl_n$ if $|Y|=n$.
\end{itemize}
\item The class of Borel subalgebras of $\sl_\infty$ appearing in (3) is precisely the class of Borel subalgebras 
parametrized by all sequential (countable) totally ordered sets $Y$.
\end{enumerate}
\eth

As pointed out in \S \ref{1.1}(a), it is important to include in our
collections of positive $c$-vectors those coming from unreachable cluster
tilting subcategories. If we do not do this, then in \thref{3} we will
obtain sets that are significantly smaller than the roots of Borel
subalgebras of $\sl_\infty$.

The paper is organized as follows. Sections \ref{cat-c}--\ref{2-formulas} contain our general results on $c$-vectors for 2-CY categories, including the 
proofs of Theorems \ref{t1} and \ref{t2}. Section \ref{backgr} contains background material on cluster categories of type $A_\infty$ and 
the Borel subalgebras of the Lie algebra $\sl_\infty$. Sections \ref{ind}--\ref{max-R+} contain our results on the $c$-vectors of the cluster categories of 
type $A_\infty$ and the indices of their objects, including the proof of \thref{3}.

\sectionnew{Homological $c$-vectors of 2-CY categories}
\label{cat-c}

\subsection{Cluster structures on triangulated $2$-CY categories}

Throughout, $\KK$ is an algebraically closed field and $\cC$ is an essentially small $\KK$-linear $\Hom$-finite triangulated category which is Krull--Schmidt and $2$-CY.  
Its suspension functor is denoted by $\Sigma$.  We recall some definitions from \cite{AusSma} and \cite{Enochs}.

\bde{homological_definitions}
Let $\cT$ be a full subcategory of $\cC$.  A {\em $\cT$-precover} of an object $c \in \cC$ is a morphism $t \xrightarrow{} c$ with $t \in \cT$, which has the following lifting property for each morphism $t' \xrightarrow{} c$ with $t' \in \cT$.  
\[
  \xymatrix 
  {
    & t \ar[d] \\
    t' \ar[r] \ar@{.>}^{\exists}[ur] & c \\
  }
\]
A $\cT$-precover $t \xrightarrow{} c$ is called a {\em $\cT$-cover} if it does not permit a non-zero direct summand of the form $t'' \xrightarrow{} 0$.  We say that $\cT$ is {\em (pre)covering} if each $c \in \cC$ has a $\cT$-(pre)cover.

The notions of {\em $\cT$-(pre)envelopes} and of $\cT$ being {\em (pre)enveloping} are defined dually.

We say that $\cT$ is {\em functorially finite} if it is precovering and preenveloping.

Note that since $\cC$ is $\Hom$-finite Krull--Schmidt, $\cT$ is covering/enveloping if and only if  it is precovering/preenveloping.

The {\em quiver} $Q_{ \cT }$ of $\cT$ was defined in \cite[sec.\ 8.1]{GabRoi}.
\ede

We assume that $\cC$ has at least one cluster tilting subcategory in the following sense.

\bde{cts}
A {\em{cluster tilting subcategory}} of $\cC$ is a functorially finite subcategory $\cT$ such that $\cT = (\Sig^{-1} \cT)^\perp = {}^\perp (\Sig \cT)$ where
\begin{align*}
\cT^\perp &:= \{ P \in \cC \mid \cC(M,P) =0 \; \; \mbox{for all} \; \; M \in \cC\},  
\\
{}^\perp \cT &:= \{ P \in \cC \mid \cC(P,M) =0 \; \; \mbox{for all} \; \; M \in \cC\}. 
\end{align*}
\ede

The following theorem recalls some important properties established in \cite[thm.\ II.1.6]{BIRS} and \cite[thm.\ 5.3]{IY}.

\bth{cluster_structures}
The category $\cC$ has a {\em weak cluster structure} in the sense that
\begin{enumerate}
\setlength\itemsep{4pt}

  
  \item  Each cluster tilting subcategory  $\cT \subseteq \cC$ can be {\em mutated} at each $t \in \Indec( \cT )$.  That is, there is a unique other isomorphism class of indecomposables $t^*$ such that $\cT^* = \add\big( ( \Indec( \cT ) \setminus t ) \cup t^* \big)$ is cluster tilting. 
  
  \item  In part (1), there are {\em exchange triangles} $t \stackrel{ \varphi }{ \longrightarrow } e \stackrel{ \varphi' }{ \longrightarrow } t^*$ and $t^* \stackrel{ \sigma }{ \longrightarrow } e' \stackrel{ \sigma' }{ \longrightarrow } t$ where $\varphi$ and $\sigma$ are $\add( \Indec( \cT ) \setminus t )$-envelopes and $\varphi'$ and $\sigma'$ are $\add( \Indec( \cT ) \setminus t )$-covers.

\end{enumerate}
If the quiver $Q_{ \cT }$ of each cluster tilting subcategory $\cT$ is without loops and $2$-cycles, then $\cC$ has a {\em cluster structure} in the sense that we also have
\begin{enumerate}
\setcounter{enumi}{2}
\setlength\itemsep{4pt}

  \item  In part (1), passing from $Q_{ \cT }$ to $Q_{ \cT^* }$ is given by {\em Fomin--Zelevinsky mutation} of $Q_{ \cT }$ at $t$, and $\dim_{ \KK } \Ext_{ \cC }( t,t^* ) = 1$.

\end{enumerate}
\eth

\subsection{Homological $g$-vectors} 
\lb{2.1}

\bde{spG}
The {\em split Grothendieck group} of an additive category will be denoted by $\Ksp(-)$. Denote by $\Indec(\cC)$ the 
isomorphism classes of indecomposable objects of $\cC$ and by $\Indec(\cT)$ those lying in $\cT$. The Krull--Schmidt assumption on $\cC$ implies that $\Ksp(\cT)$ is a free abelian group with basis $\{[t] \mid t \in \Indec(\cT) \}$:
\begin{equation}
\label{Ksp}
\Ksp(\cT) = \bigoplus_{t \in \Indec(\cT)} \Zset \cdot [t].
\end{equation}
\ede

\bde{ind-obj}
The {\em index} of an object $s \in \cC$ with respect to a cluster tilting subcategory $\cT$ is 
\[
\ind_\cT(s):= [t_0] - [t_1] \in \Ksp(\cT) \quad \mbox{where} \quad 
t_1 \to t_0 \to s \; \; \mbox{is a distinguished triangle}.
\]
\ede

The properties of indices in this general setting were investigated in \cite{DK,JP,P}. The following important 
fact was established in \cite[Theorem 2.6]{DK}.

\bth{DKthm} {\em{[Dehy--Keller]}} For every 2-CY category $\cC$ as above and two cluster tilting 
subcategories $\cT$ and $\cU$ of $\cC$, we have
\[
\Ksp(\cT) = \bigoplus_{u \in \Indec(\cU)} \Zset \cdot \ind_{\cT}(u).
\] 
\eth

\noindent
In the special case $\cU := \cT$, we recover the basis \eqref{Ksp}. We will call the cardinality $|\Indec(\cT)| = |\Indec(\cU)|$ 
rank of $\cC$.

\bre{g-vects}
We view the collection of $g$-vectors 
\[
\{ \ind_{\cT}(s) \mid \mbox{$s \in \cC$ is a rigid indecomposable object} \} 
\]
as the {\em{homological $g$-vectors of the 2-CY category $\cC$}}.  Note that
\begin{itemize}
\item $\cC$ does not necessarily have a cluster structure, and 
\item the cluster tilting subcategories of $\cC$ are not necessarily reachable from each other.
\end{itemize}
The second degree of generality is particularly important since the reachability property may fail already in finite rank, see \cite[ex.\ 4.3]{Pl}, and always fails in infinite rank. In other words the combinatorial definition of $g$-vectors \cite{FZ} 
via mutations in cluster algebras with principal coefficients reconstructs only a subset of the set of all homological $g$-vectors of a cluster category of infinite rank.
\ere

\subsection{Index maps and the dual category}
\label{2.3}
We will extensively use the homomorphism 
\begin{equation}
\label{ind-lin-map}
\ind_{\cT} : \Ksp(\cC) \to \Ksp(\cT)
\end{equation}
given by
\[
\ind_{\cT}([s]) := \ind_{\cT} (s) \quad \mbox{for} \; s \in \Indec(\cC).
\]
Consider its restriction 
\[
  \ind_{\cT}\big|_{\Ksp(\cU)} : \Ksp(\cU) \to \Ksp(\cT).
\]
The columns of the matrix of $\ind_{\cT}|_{\Ksp(\cU)}$ in the bases
\[
\{ [u] \mid u \in \Indec(\cU) \}
\quad \mbox{and} \quad 
\{ [t] \mid t \in \Indec(\cT) \}
\]
of $\Ksp(\cU)$ and $\Ksp(\cT)$ are precisely the homological $g$-vectors $\ind_{\cT}(u)$. Thus the linear map $\ind_{\cT}$ 
provides a bases-free approach to encode homological $g$-vectors. 

Denote by $\cC^{\opp}$ the opposite $2$-CY triangulated category. A full subcategory $\cT$ of $\cC$ is 
cluster tilting if and only if $\cT$ is a cluster tilting subcategory of $\cC^{\opp}$. By abuse of notation we will denote both 
categories by $\cT$. The mutation quiver of the subcategory$\cT$ of $\cC^{\opp}$ is the opposite quiver to $Q_{\cT}$. 
We will identify 
\[
\Ksp(\cC^{\opp}) \cong \Ksp(\cC).
\] 
This restricts to an identification of the split Grothendieck groups of the subcategories $\cT$ of $\cC^{\opp}$ and $\Ksp(\cC)$.
Denote by $\ol{\ind}_{\cT}(u)$ the index of $u \in \Indec(\cC^{\opp})$ and by 
\[
\ol{\ind}_{\cT} \in \Hom (\Ksp(\cC), \Ksp(\cT))
\]
the corresponding index map.
\ble{index-bar} For each cluster tilting subcategory $\cT$ of $\cC$, we have
\[
\ol{\ind}_{\cT} = - \ind_{\cT} \circ \Sigma \in \Hom (\Ksp(\cC), \Ksp(\cT)).
\]
\ele

\begin{proof}
For $c \in \cC$ it follows from \cite[proposition in sec.\ 2.1]{KR} that there is a triangle $t_1 \rightarrow t_0 \rightarrow \Sigma c$ in $\cC$ with $t_i \in \cT$ whence $-\ind_{ \cT } \circ \Sigma( c ) = [t_1] - [t_0]$.  Rolling the triangle gives a triangle $c \rightarrow t_1 \rightarrow t_0$ in $\cC$, so a triangle $t_0 \rightarrow t_1 \rightarrow c$ in the opposite category of $\cC$, whence $\overline{ \ind }_{ \cT }( c ) = [t_1] - [t_0]$.  Hence $-\ind_{ \cT } \circ \Sigma( c ) = \overline{ \ind }_{ \cT }( c )$.  
\end{proof}

\subsection{Homological $c$-vectors}
\label{2.2}
For an abelian group $A$ set
\[
A^*:= \Hom (A, \Zset).
\] 
If $\cT \subseteq \cC$ is a cluster tilting subcategory, then
\begin{equation}
\label{Ksp*}
\Ksp(\cT)^* = \prod_{t \in \Indec(\cT)} \Zset \cdot [t]^*
\end{equation}
where $[t]^* \in \Ksp(\cT)^*$ is the unique element defined by 
\[
\lcor [t]^*, [s] \rcor = \delta_{ts} \quad \mbox{for all} \quad s \in \Indec(\cT).
\]
\bde{c-vect} Let $\cT,\cU \subseteq \cC$ be cluster tilting subcategories.  For $u \in \Indec(\cU)$, define the {\em homological $c$-vector} of $(u, \cU)$ with respect to $\cT$ to be the element $c_{\cT}(u, \cU) \in \Ksp(\cT)^*$ such that 
\[
\lcor c_{\cT}(u, \cU), \ind_{\cT}(v) \rcor = \delta_{uv} \quad \mbox{for each} \quad v \in \cU. 
\]
Similarly, $\overline{c}_{\cT}(u, \cU) \in \Ksp(\cT)^*$ is the element such that 
\[
\lcor \overline{c}_{\cT}(u, \cU), \overline{\ind}_{\cT}(v) \rcor = \delta_{uv} \quad \mbox{for each} \quad v \in \cU. 
\]
\thref{DKthm} and \leref{index-bar} show that $c_{\cT}(u, \cU)$ and $\overline{c}_{\cT}(u, \cU)$ exist and are unique.
\ede
\thref{DKthm} also implies that 
\begin{equation}
\label{Ksp*2}
\Ksp(\cT)^* = \prod_{u \in \Indec(\cU)} \Zset \cdot c_{\cT}(u, \cU).
\end{equation}
Clearly, $c_{\cT}(t, \cT) = [t]^*$ for $t \in \Indec(\cT)$ and the special case of \eqref{Ksp*2} for $\cU:= \cT$ recovers 
\eqref{Ksp*}.

The next result follows from \cite{DK}, \cite{DWZ}, and \cite{NZ}.

\bth{c-vect} Let $\cT$ be a cluster tilting subcategory of 
a 2-CY category $\cC$ with a cluster structure and $\AA(Q_{\cT})$ be the corresponding 
cluster algebra. For every reachable cluster tilting subcategory $\cU$ of $\cC$ and $u \in \cU$, 
the homological $c$-vector $c_{\cT}(u, \cU) \in \Ksp(\cT)^*$ written in the basis \eqref{Ksp*} equals 
the $c$-vector of the cluster variable $x_{u} \in \AA(Q_{\cT})$ of the seed $(x_{v}, v \in \cU, Q_{\cU})$
with respect to the seed of $\AA(Q_{\cT})$ associated to $\cT$.
\eth
\begin{proof}
First observe that the $g$-vector of the cluster variable $x_v$ with respect to the seed $\Sigma$ of $\AA(Q_{\cT})$ associated to $\cT$ is given by 
\begin{equation}
\label{equ:DKDWZ}
  g_\Sigma(x_{v}) = \ind_{\cT}(v) \quad \mbox{for each} \quad v \in \Indec(\cU).
\end{equation}
This was proved in \cite[Theorem in Sect.\ 6 and end of Sect.\ 6]{DK} for a cluster algebra of finite rank satisfying the sign coherence 
conjecture for $g$-vectors of Fomin--Zelevinsky \cite[Conjectures 6.10 and 6.13]{FZ}. That conjecture was established in \cite{DWZ} in the skew-symmetric case and in \cite{GHKK} in the general case.

Denote by $\Sigma$ and $\Sigma'$ the seeds $(x_{s}, s \in \cT, Q_{\cT})$ and 
$(x_{v}, v \in \cU, Q_{\cU})$ of the cluster algebra $\AA(Q_{\cT})$ (of possibly infinite rank).
For $u \in \Indec(\cU)$, denote by $c_{\Sigma}(x_{u}, \Sigma')$ and $g_{\Sigma}(x_{u})$ 
the $c$-vector of the pair $(x_{u}, \Sigma')$ with respect to the seed $\Sigma$ and the 
$g$-vector of the cluster variable $x_{u}$ with respect to $\Sigma$. Both are viewed 
as elements of $\Ksp(\cT)$ in the basis \eqref{Ksp}. 

We embed $\Ksp(\cT) \hra \Ksp(\cT)^*$ by sending $[t] \mt [t]^*$ for $t \in \Indec(\cT)$. Eq. (1.11) of 
\cite[Theorem 1.2]{NZ} implies that, under this embedding,
\[
\lcor c_\Sigma(x_{u}, \Sigma'), g_\Sigma(x_{v}) \rcor = \delta_{uv} \quad \mbox{for all} \quad v \in \cU. 
\]
Combining with Equation \eqref{equ:DKDWZ} shows $c_\Sigma(x_{u}, \Sigma') = c_{\cT}(u, \cU)$, which completes the proof of the theorem.
\end{proof}

\bre{g-vects2} \hfill
\begin{enumerate}
\item It follows from \thref{c-vect} that for all reachable (form $\cT$) 
cluster tilting subcategories $\cU$ of $\cC$, the $c$-vectors $c_{\cT}(u, \cU) \in \Ksp(\cT)^*$ 
are finite combinations of $\{[t]^* \mid t \in \Indec(\cU)\}$. We will see in Sect. \ref{icvect} that this is not the case for
the unreachable cluster tilting subcategories $\cU$ of $\cC$.
\item Similarly to the case of $g$-vectors, the combinatorial definition of $c$-vectors \cite{FZ} via mutations 
in cluster algebras with principal coefficients reconstructs only a subset of the set of all homological $c$-vectors 
of a cluster category of infinite rank.
\item The $c$-vector $c_{\cT}(u, \cU)$ depends on the pair $(u, \cU)$ and not only on 
the indecomposable object $u$. This follows from from \thref{c-vect} and the analogous 
known fact for the $c$-vectors obtained via mutation \cite{FZ}.  
\end{enumerate}
\ere

\ble{c_versus_c_bar}
If $\cT,\cU \subseteq \cC$ are cluster tilting subcategories and $u \in \cU$ is indecomposable, then $c_{ \cT }( \Sigma u,\Sigma\cU ) = - \overline{ c }_{ \cT }( u,\cU )$.
\ele

\begin{proof}
If $v \in \cU$ is indecomposable, then
\[
  \langle \overline{ c }_{ \cT }( u,\cU ) , \overline{ \ind }_{ \cT }( v ) \rangle = \delta_{ uv }
\]
by definition, while
\[
  \langle c_{ \cT }( \Sigma u,\Sigma \cU ),\overline{ \ind }_{ \cT }( v ) \rangle =
  - \langle c_{ \cT }( \Sigma u,\Sigma \cU ),\ind_{ \cT }( \Sigma v ) \rangle = - \delta_{ uv }
\]
where the first equality is by Lemma \ref{lindex-bar} and the second is by definition.  This proves the lemma since $\overline{ \ind }_{ \cT }( v )$ ranges through a basis of $\Ksp(\cT)$.
\end{proof}

%
%
%
%
%
%
%
%

\sectionnew{Decomposing the set of $c$-vectors}

\subsection{Sign coherence of $c$-vectors}

\bde{sign_coherence}
For $a,b \in \Ksp(\cT)^*$ we write $a \leqslant b$ if
\[
  \langle a,[t] \rangle \leqslant \langle b,[t] \rangle
  \quad \mbox{for each} \; t \in \Indec( \cT ).
\]
The inequality $a \geqslant b$ is defined analogously.

An element $c \in \Ksp(\cT)^*$ is called {\em sign coherent} if it satisfies $c \leqslant 0$ or $c \geqslant 0$, that is, if each non-zero $\langle c,[t] \rangle$ with $t \in \Indec( \cT )$ has the same sign.

Note that $c = 0$ is sign coherent.
\ede

\bde{C_sets}
Let $\cT \subseteq \cC$ be a cluster tilting subcategory. 
\begin{enumerate}
\setlength\itemsep{4pt}

  \item  $C_{\cT}(\cC) = \{\, c_{ \cT }( u,\cU ) \,|\, \mbox{$\cU \subseteq \cC$ cluster tilting and $u \in \Indec( \cU )$} \,\}$.

  \item  $C^+_{\cT}(\cC) = \{\, c \in C_{\cT}(\cC) \,|\, c \geqslant 0 \,\}$.
  
  \item  $C^-_{\cT}(\cC) = \{\, c \in C_{\cT}(\cC) \,|\, c \leqslant 0 \,\}$.
  
\end{enumerate}
The sets $\overline{C}_{\cT}(\cC)$ and $\overline{C}^{\pm}_{\cT}(\cC)$ are defined analogously with $\overline{c}_{ \cT }( u,\cU )$ instead of $c_{ \cT }( u,\cU )$.
\ede

\bre{non-zero}
Note that $0 \notin C_{\cT}(\cC)$ because of \thref{DKthm} and the definition of $c$-vectors.
\ere

\bpr{sign_coherence}
\begin{enumerate}
\setlength\itemsep{4pt}

  \item  If $\cT,\cU \subseteq \cC$ are cluster tilting subcategories and $u \in \cU$ is indecomposable, then $c_{ \cT }( u,\cU )$ and $\overline{ c }_{ \cT }( u,\cU )$ are sign coherent.

  \item  If $\cT \subseteq \cC$ is a cluster tilting subcategory then $C_{ \cT }( \cC ) = C_{ \cT }^+( \cC ) \sqcup C_{ \cT }^-( \cC )$ and $\overline{ C }_{ \cT }( \cC ) = \overline{ C }_{ \cT }^+( \cC ) \sqcup \overline{ C }_{ \cT }^-( \cC )$.
  
\end{enumerate}

\epr

\begin{proof}
Part (2) is a reformulation of part (1).  By Lemma \ref{lc_versus_c_bar} it is enough to prove part (1) for $\overline{ c }_{ \cT }( u,\cU )$.  Let $s,t \in \cT$ be indecomposable.  Since $s \oplus t$ is rigid, \cite[sec.\ 2.4]{DK} says that if $[u]$ appears with non-zero coefficients in $\ind_{ \cU }( s )$ and $\ind_{ \cU }( t )$ then the coefficients have the same sign.  That is, $[u]^{ \ast } \circ \ind_{ \cU }( s ) = \langle \overline{ c }_{ \cT }( u,\cU ),s \rangle$ and $[u]^{ \ast } \circ \ind_{ \cU }( t ) = \langle \overline{ c }_{ \cT }( u,\cU ),t \rangle$ have the same sign as desired.
\end{proof}

The next proposition completes our results on decompositions of the sets of $c$-vectors. Its proof relies on \thref{2}(1)
which is independently proved in the next section.

\bpr{positive_and_negative_c-vectors}
Assume that the quiver $Q_{ \cT }$ of each cluster tilting subcategory $\cT$ is without loops and $2$-cycles.

Let $\cT \subseteq \cC$ be a cluster tilting subcategory.
\begin{enumerate}
\setlength\itemsep{4pt}

  \item  $C_{ \cT }^+( \cC ) = -C^-_{ \cT }( \cC ) = \overline{ C }^+_{ \cT }( \cC ) = - \overline{ C }^-_{ \cT }( \cC )$.

  \item  $C_{ \cT }( \cC ) = \overline{ C }_{ \cT }( \cC )$.

\end{enumerate}
\epr

\begin{proof}
(1)  Lemma \ref{lc_versus_c_bar} implies $C^+_{ \cT }( \cC ) = -\overline{ C }^-_{ \cT }( \cU )$ and $- C^-_{ \cT }( \cC ) = \overline{ C }^+_{ \cT }( \cC )$.  It is hence sufficient to prove $\overline{ C }^+_{ \cT }( \cC ) = -\overline{ C }^-_{ \cT }( \cC )$.  We only show $\subseteq$ since $\supseteq$ follows by a symmetric argument.

Consider an element of $\overline{ C }^+_{ \cT }( \cC )$.  It has the form $\overline{ c }_{ \cT }( u,\cU ) = [u]^{ \ast } \circ \ind_{ \cU }\big|_{ \Ksp( \cT ) }$ where $\cU \subseteq \cC$ is a cluster tilting subcategory and $u \in \cU$ is indecomposable.  We must have
\begin{equation}
\label{equ:positive_and_negative_c-vectors}
  [u]^{ \ast } \circ \ind_{ \cU }\big( [t] \big) \geqslant 0	
\end{equation}
for each $t \in \cT$.

We can mutate $\cU$ at $u$ and let $\cU^{ \ast }$, $u^{ \ast }$ be the mutated category and object.  Let $u^{ \ast } \rightarrow e \rightarrow u$ be the corresponding exchange triangle and consider the homomorphism
$\phi_+ : \Ksp( \cU ) \rightarrow \Ksp( \cU^{ \ast } )$ given by
\[
  \phi_+\big( [ \widetilde{u} ] \big)
  =
  \left\{
    \begin{array}{cl}
      {[e]} - [ u^{ \ast } ] & \mbox{ if $\widetilde{ u } \cong u$, } \\[2mm]
      [ \widetilde{u} ] & \mbox{ if $\widetilde{ u } \in \cU$ is an indecomposable not isomorphic to $u$, }
    \end{array}
  \right.
\]
which was introduced in \cite[sec.\ 3]{DK}.  It follows from \cite[thm.\ 3.1]{DK} and Equation \eqref{equ:positive_and_negative_c-vectors} that $\ind_{ \cU^{ \astsmall } }\big( [t] \big) = \phi^+ \circ \ind_{ \cU }\big( [t] \big)$ for each $t \in \cT$.    This implies (a) in the following computation:
\begin{align*}
  \langle \overline{ c }_{ \cT }( u^{ \ast },\cU^{ \ast } ),[t] \rangle
  & = [ u^{ \ast } ]^{ \ast } \circ \ind_{ \cU^{ \astsmall } }\big( [t] \big) \\
  & \stackrel{ \rm (a) }{ = }
    [ u^{ \ast } ]^{ \ast } \circ \phi^+ \circ \ind_{ \cU }\big( [t] \big) \\
  & \stackrel{ \rm (b) }{ = }
    - [ u ]^{ \ast } \circ \ind_{ \cU }\big( [t] \big) \\
  & = - \langle \overline{ c }_{ \cT }( u,\cU ),[t] \rangle,
\end{align*}
where (b) follows from the formula for $\phi^+$.  Hence
\[
  \overline{ c }_{ \cT }( u,\cU )
  = - \overline{ c }_{ \cT }( u^{ \ast },\cU^{ \ast } )
  \in -\overline{ C }^-_{ \cT }( \cC )
\]
as desired.

(2)  Immediate by part (1) and Proposition \ref{psign_coherence}(2).
\end{proof}


\subsection{A decomposition of the set of positive $c$-vectors}
\label{sec:decom}
\bde{R}
Let $P$ be a partially ordered set. A non-empty subset $I \subset P$ is called an ideal if the following two conditions hold:
\begin{enumerate}
\setlength\itemsep{4pt}

  \item  If $c,c' \in P$ have $c \geqslant c'$, then $c \in I \Rightarrow c' \in I$.
  
  \item If $c_1,\ldots,c_n \in I$ are given then there exists a $c \in I$ with $c \geqslant c_j$ for each $j$.

\end{enumerate}  
By a maximal ideal of $P$ we will mean an ideal that is maximal with respect to inclusion, but is not necessarily proper.
\ede

For each cluster tilting subcategory $\cT \subseteq \cC$, $C^+_{ \cT }( \cC )$ has a natural poset structure by \deref{sign_coherence}.
$C^+_{ \cT }( \cC )$ is the union of its maximal ideals with respect to this poset structure.

\sectionnew{Two effective formulae for $c$-vectors}
\label{2-formulas}

This section proves Theorems \ref{t1} and \ref{t2} from the introduction.

\subsection{Proof of Theorem \ref{t1}}
By Lemma \ref{lindex-bar} it is enough to prove that there are inverse isomorphisms
\[
  \xymatrix {
  \Ksp( \cT )
    \ar[rrrr]<1ex>^-{ -\ind_{ \cU } \circ \Sigma \big|_{ \Ksp( \cT ) } } & & & &
    \Ksp( \cU ).
    \ar[llll]<1ex>^-{ \ind_{ \cT } \big|_{ \Ksp( \cU ) } }
            }
\]
For $u \in \cU$, it was shown in \cite[proposition in sec.\ 2.1]{KR} that there is a triangle $t_1 \rightarrow t_0 \rightarrow u \rightarrow \Sigma t_1$ with $t_i \in \cT$.  By definition, $\ind_{ \cT }\big( [u] \big) = [t_0] - [t_1]$ whence
\[
  \Sigma \circ \ind_{ \cT }\big( [u] \big) = [\Sigma t_0 ] - [\Sigma t_1].
\]
Turning the triangle gives $u \rightarrow \Sigma t_1 \rightarrow \Sigma t_0 \stackrel{ \delta }{ \rightarrow } \Sigma u$.  The morphism $\Hom_{ \cC }( -,\Sigma t_0 )\big|_{ \cU } \rightarrow \Hom_{ \cC }( -,\Sigma u )\big|_{ \cU }$ induced by $\delta$ is zero because its target is zero, so \cite[prop.\ 2.2]{P} gives
\[
  \ind_{ \cU }\big( [u] \big) =
  \ind_{ \cU }\big( [\Sigma t_1] \big)
  - \ind_{ \cU }\big( [\Sigma t_0] \big).
\]
The displayed formulae give (a) and (b) in the following computation.
\begin{align*}
  [u] & = \ind_{ \cU }\big( [u] \big) \\
      & \stackrel{ \rm (a) }{ = } \ind_{ \cU }\big( [\Sigma t_1] \big) - \ind_{ \cU }\big( [\Sigma t_0] \big) \\
      & = - \ind_{ \cU }\big( [ \Sigma t_0 ] - [ \Sigma t_1 ] \big) \\
      & \stackrel{ \rm (b) }{ = } - \ind_{ \cU } \circ \Sigma \circ \ind_{ \cT }\big( [u] \big)
\end{align*}
This shows that one of the compositions of the homomorphisms in the theorem is the identity.  The theorem follows because the second homomorphism, $\ind_{ \cT }\big|_{ \Ksp( \cU ) }$, is an isomorphism by \cite[thm.\ 2.4]{DK}.

\subsection{Proof of Theorem \ref{t2}(1)}
We can compute for $u,v \in \Indec( \cU )$:
\begin{align*}
  \langle [u]^* \circ \overline{ \ind }_{ \cU } \big|_{ \Ksp( \cT ) },\ind_{ \cT }( v ) \rangle
  & = [u]^* \circ \overline{ \ind }_{ \cU } \big|_{ \Ksp( \cT ) } \circ \ind_{ \cT }( v ) \\
  & \stackrel{ \rm (a) }{ = } [u]^*\big( [v] \big) \\
  & = \delta_{ uv },
\end{align*}
where (a) is by Theorem \ref{t1}.  By Definition \ref{dc-vect} this implies the formula in Theorem \ref{t2}(1).

\subsection{Proof of Theorem \ref{t2}(2)}

\bre{c-vectors_are_dimension_vectors}
Let $\cT,\cU \subseteq \cC$ be cluster tilting subcategories and let $t \in \cT$, $u \in \cU$ be given with $u$ indecomposable.  Let $\cU^{ \ast }$ and $u^{ \ast }$ be defined by mutation of $\cU$ at $u$.  There is an additive subcategory $\widetilde{ \cU }$ such that $\cU = \add( \widetilde{ \cU } \cup \{ u \} )$ and $\cU^{ \ast } = \add( \widetilde{ \cU } \cup \{ u^{ \ast } \} )$.  Let
\begin{equation}
\label{equ:exchange_triangle}	
  u^{ \ast } \rightarrow e \rightarrow u \stackrel{ \delta }{ \rightarrow } \Sigma u^{ \ast }
\end{equation}
be an exchange triangle with respect to $\cU$.  There is a morphism
\bgroup
\[
  \varphi :
    \Hom_{ \cC }( -,u )\big|_{ \cT }
    \rightarrow
    \Hom_{ \cC }( -,\Sigma u^{ \ast } )\big|_{ \cT }
\]
\egroup
induced by $\delta$.  Let
\begin{equation}
\label{equ:u_t_triangle}
	t \stackrel{ \alpha }{ \rightarrow } u^0 \rightarrow u^1
\end{equation}
be a triangle in $\cC$ with $u^i \in \cU$, see \cite[proposition in sec.\ 2.1]{KR}.  There is an induced map
\bgroup
\[
  \alpha^{ \ast } :
    \Hom_{ \cC }( u^0,\Sigma u^{ \ast } )
    \rightarrow
    \Hom_{ \cC }( t,\Sigma u^{ \ast } ).
\]
\egroup
\ere

\ble{identical_images}
The induced maps $\varphi_t$ and $\alpha^{ \ast }$ have identical images.
\ele

\begin{proof}
The triangles \eqref{equ:exchange_triangle} and \eqref{equ:u_t_triangle} give the following commutative diagram with an exact row and an exact column.
\bgroup
\[
  \vcenter{
  \xymatrix @+1.0pc 
                       {
    & \Hom_{ \cC }( t,u ) \ar[r] \ar_{ \varphi_t }[d] & \Hom_{ \cC }( \Sigma^{ -1 }u^1,u ) \ar[d] \\
    \Hom_{ \cC }( u^0,\Sigma u^{ \ast } ) \ar^{ \alpha^{ \astsmall } }[r] \ar[d] & \Hom_{ \cC }( t,\Sigma u^{ \ast } ) \ar[r] \ar[d] & \Hom_{ \cC }( \Sigma^{ -1 }u^1,\Sigma u^{ \ast } ) \\
    \Hom_{ \cC }( u^0,\Sigma e ) \ar[r] & \Hom_{ \cC }( t,\Sigma e ) \\
                        }
          }
\]
\egroup
The northeast and southwest objects, $\Hom_{ \cC }( \Sigma^{ -1 }u^1,u )$ and
$\Hom_{ \cC }( u^0,\Sigma e )$, are zero because $u,u_0,u_1,e \in \cU$, and this implies the lemma.
\end{proof}

\ble{c-vectors_are_dimensions_of_images}
Assume that the quiver of each cluster tilting subcategory is without loops and $2$-cycles.  If $c_{ \cT }( u,\cU ) \geqslant 0$ then 
\[
  \langle c_{ \cT }( u,\cU ),t \rangle
  = \dim_{ \KK }\!
    (
      \Image \alpha^{ \ast }
    ).
\]
\ele

\begin{proof}
We can rewrite the triangle \eqref{equ:u_t_triangle} in $\cC$ as
\begin{equation}
\label{equ:refined_triangle}
  t
  \stackrel{ \alpha }{ \rightarrow }
  \widetilde{ u }^0 \oplus u^m
  \rightarrow 
  \widetilde{ u }^1 \oplus u^n
\end{equation}
with $\widetilde{ u }^i \in \widetilde{ \cU }$.  This can be interpreted as a triangle $\widetilde{ u }^1 \oplus u^n \rightarrow \widetilde{ u }^0 \oplus u^m \stackrel{ \alpha }{ \rightarrow } t$ in  the opposite category of $\cC$, so $[u]^* \circ \overline{ \ind }_{ \cU }( t ) = m - n$.  Hence Theorem \ref{t2}(1) implies
\[
  \langle c_{ \cT }( u,\cU ),t \rangle = m - n.
\]
It follows from \cite[prop.\ 2.1]{DK} that at most one of $m$ and $n$ is non-zero, and it must be $m$ because $c_{ \cT }( u,\cU ) \geqslant 0$.  Hence $n = 0$ and
\begin{equation}
\label{equ:c_bar_is_m}
  \langle c_{ \cT }( u,\cU ),t \rangle = m.
\end{equation}

On the other hand, acting on \eqref{equ:refined_triangle} with $\Hom_{ \cC }( -,\Sigma u^{ \ast } )$ gives an exact sequence
\[
  \Hom_{ \cC }( \widetilde{ u }^1 \oplus u^n,\Sigma u^{ \ast } )
  \longrightarrow
  \Hom_{ \cC }( \widetilde{ u }^0 \oplus u^m,\Sigma u^{ \ast } )
  \stackrel{ \alpha^{ \astsmall } }{ \longrightarrow }
  \Hom_{ \cC }( t,\Sigma u^{ \ast } ).
\]
Note that $\Hom_{ \cC }( \widetilde{ u }^i,\Sigma u^{ \ast } ) = 0$ because $\widetilde{ u }^i,u^{ \ast }, \in \cU^{ \ast }$.  Moreover, we proved $n = 0$ and have $\Hom_{ \cC }( u,\Sigma u^{ \ast } ) \cong \KK$ by \thref{cluster_structures}(3).  Hence the exact sequence is isomorphic to
\[
  0
  \longrightarrow
  \KK^m
  \stackrel{ \alpha^{ \astsmall } }{ \longrightarrow }
  \Hom_{ \cC }( t,\Sigma u^{ \ast } ),
\]
so
\[
  \dim_{ \KK } ( \Image \alpha^{ \ast } ) = m.
\]
Combining with Equation \eqref{equ:c_bar_is_m} concludes the proof.
\end{proof}

We can now finally prove Theorem \ref{t2}(2): Assume that the quiver of each cluster tilting subcategory is without loops and $2$-cycles.

If $c_{ \cT }( u,\cU ) \geqslant 0$ then
\[
  \langle \dim_{ \cT } ( \Image \varphi ),t \rangle
  = \dim_{ \KK }( \Image \varphi_t )
  = \dim_{ \KK } ( \Image \alpha^{ \ast } )
  = \langle c_{ \cT }( u,\cU ),t \rangle,
\]
where the equalities are, respectively, by definition, by Lemma \ref{lidentical_images}, and by Lemma \ref{lc-vectors_are_dimensions_of_images}.  Since $c_{ \cT }( u,\cU ) \neq 0$ it follows that $\Image \varphi \neq 0$ whence $\delta \not\in [\Sigma \cT]$.

On the other hand, if $c_{ \cT }( u,\cU ) \not\geqslant 0$ then $c_{ \cT }( u,\cU ) \leqslant 0$ by Proposition \ref{psign_coherence} and $c_{ \cT }( u^{ \ast },\cU^{ \ast } ) \geqslant 0$ by the proof of Proposition \ref{ppositive_and_negative_c-vectors}.  The previous paragraph shows $\delta^{ \ast } \not\in [\Sigma \cT]$ where $u \xrightarrow{} e' \xrightarrow{} u^{ \ast } \xrightarrow{ \delta^{ \ast } } \Sigma u$ is an exchange triangle.  By \cite[prop.\ 4.3]{P} this implies $\delta \in [\Sigma \cT]$.

\sectionnew{Background on cluster categories of type $A_\infty$ and the Lie algebra $\sl_\infty$}
\lb{backgr}
\subsection{Cyclically ordered sets and triangulations}
Let $Z$ be a {\em{cyclically ordered set}} with $|Z| \geqslant 4$, realized as a subset of $S^1$. We will assume that it satisfies the 
Igusa--Todorov conditions:
\begin{enumerate}
\item {\em{No point of $Z$ is a limit point of $Z$}}. 
\item {\em{Each limit point of $Z$ is a limit of two sequences of points in $Z$ approaching it from both sides.}}
\end{enumerate}
These conditions imply that each $z \in Z$ has an immediate successor and an immediate predecessor, which will be denoted by 
$z^+$ and $z^-$, respectively. These are the nearest counterclockwise and clockwise points in $Z$. The elements of $Z$ will be
called {\em vertices}.  An {\em arc} of $Z$ is a set $\{ z_1,z_2 \}$ with $z_1 \neq z_2$ in $Z$.  An arc will be called an {\em edge} if the $z_i$ 
are neighbouring vertices; otherwise it will be called a {\em diagonal}.  A {\em triangle of $Z$} is a subset $\{ z_1,z_2,z_3 \} \subset Z$ with three elements such that each of $\{ z_1,z_2 \}$, $\{ z_1,z_3 \}$, $\{ z_2,z_3 \}$ is an edge or a diagonal of $Z$.  

For $x,y \in S^1$, denote by $\llbracket x,y \rrbracket$ the subset of $S^1$ 
traced from $x$ to $y$ in the counterclockwise direction.  It can be viewed as a closed interval of $\BR$, whence each subset has a supremum and an infimum.  The supremum and infimum of $\llbracket x,y \rrbracket \cap U$ for $U \subseteq S^1$ will be denoted by
\[
\sup_{ \llbracket z_1,z_2 \rrbracket }U \; \; \mbox{and} \; \; \inf_{ \llbracket z_1,z_2 \rrbracket }U.
\]
These belong to $\overline{U}$.  We will denote by $e_j \nearrow p$ and $e_j \searrow p$ the convergence in $\llbracket z_1,z_2 \rrbracket$ from below and above, for $e_j$ and $p$ in $\llbracket z_1,z_2 \rrbracket$. In terms of the circle $S^1$, this corresponds to convergence in the counterclockwise and clockwise direction, respectively.

\subsection{The cluster categories $\cC( Z )$ of type $A_\infty$}
\label{3.1}
Starting from a set $Z$ satisfying the conditions (1) and (2), Igusa and Todorov constructed in \cite[thm.\ 2.4.1]{IT2} a cluster category $\cC( Z )$ of type $A_{ \infty }$.  It is a $\KK$-linear Hom-finite Krull--Schmidt triangulated category which is 2-CY.  It is defined as follows. Set $R:= \KK[[t]]$. Let $R Z$ be the $R$-category with objects indexed by $Z$ and morphism spaces
\[
\PP(z_1, z_2) := R f_{z_2 z_1} 
\]
such that $f_{z_3 z_2} f_{z_2 z_1} = t^c f_{z_3 z_1}$, where $c =1$ if $z_1 \leqslant z_3 < z_2$ and $c=0$ otherwise. 
The map $Z \to Z$, given by $z \mt z^+$, is an isomorphism of cyclically ordered sets. It induces an autoequivalence of 
$\add R Z$ which will be denoted by $\eta$. Denote by $\MM \FF(Z)$ the exact category of pairs $(P, d)$, where $P \in \add RZ$ 
and $d \colon P \to P$ is a morphism that satisfies $d^2 = (\cdot t)$ and factors through $\eta_P$.  Igusa and Todorov 
proved that $\MM \FF(Z)$ is a Frobenius category; $\cC(Z)$ is defined to be the corresponding stable category. The latter 
has finite dimensional Hom spaces over $\KK$. The indecomposable objects 
of $\cC(Z)$ are indexed by the diagonals of $Z$. For two non-neighbouring vertices $z_1, z_2 \in Z$, the associated indecomposable object 
of $\cC(Z)$ will be denoted by $\{ z_1, z_2 \}$. The suspension functor of $\cC(Z)$ is given by $\Sig \{z_1, z_2\}= \{z_1^-, z_2^-\}$. 

We note that the categories $\cC(Z)$ and $\MM \FF(Z)$ were denoted by $\cC_\phi(Z)$ and $\MM \FF_\phi(Z)$ in \cite{IT2} 
where $\phi \in \Aut(Z)$ is given by $\phi (z) := z^+$. The indecomposable objects of the category $\cC(Z)$ were denoted by $E(z_1,z_2)$ in \cite{IT2}.
\subsection{Cluster structure of the categories $\cC( Z )$} 
\label{3.2}
A maximal collection $T$ of pairwise noncrossing diagonals of $Z$ will be called a {\em{weak triangulation of $Z$}} (more precisely, it is a weak triangulation of the $\infty$-gon with vertices in $Z$). 

The limit point $y$ of $Z$ is said to be {\em{surrounded by a fountain of $T$}} if $T$ contains diagonals $\{z_n, z\}$ and  $\{z'_n, z\}$ $(n \in \Nset)$ for some $z \in Z$ such that the sequences $z_n$ and $z'_n$ approach $y$ from the two different sides.  The limit point $y$ of $Z$ is said to be {\em{surrounded by a leapfrog of $T$}} if $T$ contains diagonals $\{z_n, z_{n+1}\}$ $(n \in \Nset)$ such that the sequences $z_{2n}$ and $z_{2n+1}$ approach $y$ from the two different sides. 

The {\em{cluster tilting subcategories}} of $\cC( Z )$ were classified in \cite[thm.\ 0.5]{GHJ} and \cite[thm. 7.17]{SvR}, 
where it was proved that they are in bijection with the weak triangulations $T$ of $Z$ that satisfy the following condition:

\medskip
(3) {\em{Each limit point $y$ of $Z$ is surrounded either by a fountain or a leapfrog of $T$.}}
\medskip

By a {\em{triangulation of $Z$}} we will mean a weak triangulation which satisfies this property.  The cluster tilting subcategory $\cT(T) \subset \cC( Z )$ corresponding to a triangulation $T$ is the additive envelope of the indecomposable objects $\{z_1, z_2\}$ where  $\{z_1, z_2\}$ runs over all diagonals in the triangulation $T$.   

It was proved in \cite[thm.\ 0.6]{GHJ} that the quiver $Q_{ \cT(T) }$  of the category $\cT(T)$ is without loops and $2$-cycles, and that hence $\cC( Z )$ has a cluster structure as explained in \thref{cluster_structures}. It is isomorphic to the dual quiver of the triangulation $T$.
	
\begin{figure}
  \centering
    \begin{tikzpicture}[scale=2.5]
      \draw (0,0) circle (1cm);

      \draw[very thick,green] ([shift=(45:0.95cm)]0,0) arc (45:200:0.95cm);
      \draw[very thick,green] ([shift=(230:0.95cm)]0,0) arc (230:375:0.95cm);
      \draw[very thick,blue] ([shift=(45:0.90cm)]0,0) arc (45:170:0.90cm);
      \draw[very thick,blue] ([shift=(230:0.90cm)]0,0) arc (230:340:0.90cm);

      \draw (45:0.97cm) -- (45:1.03cm);
      \draw (45:1.13cm) node{$x_0$};
      \draw (85:0.97cm) -- (85:1.03cm);
      \draw (85:1.13cm) node{$s_0$};
      \draw (170:0.97cm) -- (170:1.03cm);
      \draw (170:1.13cm) node{$y_0$};
      \draw (200:0.97cm) -- (200:1.03cm);
      \draw (200:1.20cm) node{$x_1^{--}$};
      \draw (215:0.97cm) -- (215:1.03cm);
      \draw (215:1.15cm) node{$x_1^-$};      
      \draw (230:0.97cm) -- (230:1.03cm);
      \draw (230:1.13cm) node{$x_1$};
      \draw (300:0.97cm) -- (300:1.03cm);
      \draw (300:1.13cm) node{$s_1$};
      \draw (340:0.97cm) -- (340:1.03cm);
      \draw (340:1.13cm) node{$y_1$};
      \draw (375:0.97cm) -- (375:1.03cm);
      \draw (375:1.20cm) node{$x_0^{--}$};
      \draw (390:0.97cm) -- (390:1.03cm);
      \draw (390:1.17cm) node{$x_0^-$};

    \end{tikzpicture} 
  \caption{Let $x = \{ x_0,x_1 \} \in \cC( Z )$ be given.  There is a non-zero morphism $x \rightarrow y$ when $y = \{ y_0,y_1 \}$ has a vertex in each green interval.  It factorises $x \rightarrow s \rightarrow y$ when $s = \{ s_0,s_1 \}$ has a vertex in each blue interval.}
\label{fig:Igusa-Todorov_morphisms}
\end{figure}
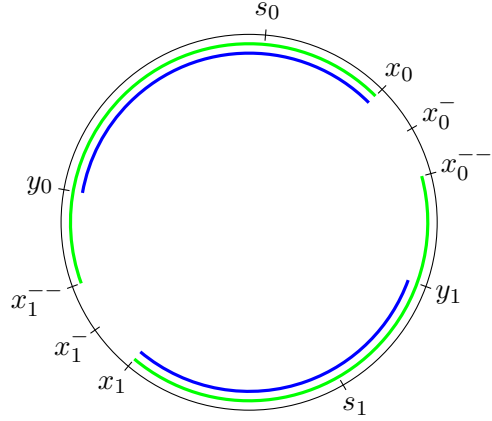

The split Grothendieck group of the cluster tilting subcategory $\cT(T)$ associated to the triangulation $T$ is 
\begin{equation}
\label{K-inft}
\Ksp(\cT(T)) := \bigoplus_{ \{z_1,z_2\} \in T} \Zset \cdot [ \{z_1, z_2\} ]
\end{equation}
and its dual is 
\begin{equation}
\label{k-inft*}
\Ksp(\cT(T))^* := \prod_{ \{z_1,z_2\} \in T} \Zset \cdot [ \{z_1, z_2\} ].
\end{equation}
\bre{Igusa-Todorov_morphisms}
If $x,y$ are indecomposable objects of $\cC( T )$ then
\[
  \dim_{ \KK }\Ext^1_{ \cC( Z ) }( x,y )
  =
  \left\{
    \begin{array}{cl}
      $1$ & \mbox{ if the diagonals of $x$ and $y$ cross, } \\
      $0$ & \mbox{ if not. }
    \end{array}
  \right.
\]
This can be refined.  Let $x,s,y$ be indecomposable objects in $\cC( Z )$.
\begin{enumerate}

  \item  There is a non-zero morphism $x \rightarrow y$ if and only if $x = \{ x_0,x_1 \}$ and $y = \{ y_0,y_1 \}$ with $x_0 \leqslant y_0 \leqslant x_1^{--} < x_1 \leqslant y_1 \leqslant x_0^{--}$.

  \item  A non-zero morphism $x \rightarrow y$ factorises $x \rightarrow s \rightarrow y$ if and only if $x = \{ x_0,x_1 \}$, $s = \{ s_0,s_1 \}$, $y = \{ y_0,y_1 \}$, with $x_0 \leqslant s_0 \leqslant y_0 \leqslant x_1^{--} < x_1 \leqslant s_1 \leqslant y_1 \leqslant x_0^{--}$.
  
\end{enumerate}
See Figure \ref{fig:Igusa-Todorov_morphisms}.
\ere

Finally, we recall that a vertex $e \in Z$ is called an {\em ear} of a triangulation $T$, if $\{e^-,e^+\}$ is a diagonal in $T$.
\subsection{Borel subalgebras of $\sl_\infty$} 
\label{3.3}
The Lie algebra $\sl_\infty:=\sl_\infty(\KK)$ is the Lie algebra of $\Zset \times \Zset$ traceless matrices with only finitely many nonzero entries
in a field $\KK$. It can be equivalently 
defined as the Lie algebra of $\Nset \times \Nset$ traceless matrices with only finitely many nonzero entries. We will use the first presentation. 
The standard {\em{Cartan subalgebra}} of $\sl_\infty$ is the diagonal subalgebra
\[
\hh : = \{ \diag(\ldots, c_{-1}, c_0, c_1, \ldots) \mid c_i \in \KK, c_i=0 \; \mbox{for} \; |i| \gg 0, \sum c_i = 0 \}. 
\]
The {\em{roots}} of $\sl_\infty$ are 
\[
\De_{\sl_\infty} := \{ \ep_i - \ep_j \mid i \neq j \in \Zset \}
\]
where $\ep_i \in \hh^*$ are the functionals
\[
\ep_i ( \diag(\ldots, c_{-1}, c_0, c_1, \ldots) ):= c_i.
\]
The corresponding {\em{root spaces}} are 
\[
\sl_\infty^{\ep_i - \ep_j} = \KK E_{ij}.
\]
For every bijection $\sigma : \Zset \to \Zset$, the map
\[
E_{ij} \mt E_{\sigma(i) \sigma(j)}
\] 
is an automorphism of $\sl_\infty$, preserving $\hh$.
For a countable set $Y$, denote
\begin{equation}
\label{Z-Y}
\Zset^{(Y)} := \bigoplus_{y \in Y} \Zset \ep_y 
\quad \mbox{and} \quad 
\Zset^{(Y)}_0 = \Span \{ \ep_y - \ep_{y'} \mid y \neq y' \in Y \} \subset \Zset^{(Y)}.
\end{equation}
We will identify
\[
\De_{\sl_\infty} \cong \{ \ep_y - \ep_{y'} \mid y \neq y' \in Y \}
\]
using any bijection $Y \cong \Zset$.

A {\em{triangular decomposition}} \cite{DP1,DP2} of $\De_{\sl_\infty}$ is a decomposition of the form
\begin{equation}
\label{triangular}
\De_{\sl_\infty} := \De^+_{\sl_\infty} \sqcup ( - \De^+_{\sl_\infty})
\end{equation}
for a subset $\De^+_{\sl_\infty} \subset \De_{\sl_\infty}$ such that the cone $\Rset_+ \De^+_{\sl_\infty}$ does not contain a line. 
A {\em{Borel subalgebra}} of $\sl_\infty$ is a subalgebra of the form \cite{DP1,DP2} 
\[
\b(\De^+_{\sl_\infty}):= \hh \oplus \left( \oplus_{\al \in \De^+_{\sl_\infty}} \sl_\infty^\al \right)
\]
for a given triangular decomposition \eqref{triangular}.

Unlike the finite dimensional case, the Borel subalgebras of $\sl_\infty$ are generally non-isomorphic to each other. The set of all triangular 
decompositions \eqref{triangular} is in bijection with the countable totally ordered sets \cite[prop. 2]{DP1}. 
A bijection of $\Zset$ sends a triangular decomposition of $\Delta_{\sl_\infty}$ to a triangular decomposition of $\Delta_{\sl_\infty}$.
Two triangular decompositions are considered equivalent if they are obtained from each other in this way.
Two countable totally ordered sets are 
considered equivalent if they are isomorphic or reverse-isomorphic.
The triangular decomposition corresponding to a countable totally ordered set $Y$ corresponds to the set 
\begin{equation}
\label{De-Y}
\De^+_{\sl_\infty}(Y) := \{ \ep_{y} - \ep_{y'} \mid y > y' \in Y \}.
\end{equation}
\sectionnew{A formula for indices in $\cC( Z )$}
\lb{ind}
In this section we prove a formula for the index of an indecomposable object of a cluster category $\cC( Z )$ of type $A_{ \infty }$ with respect to an arbitrary 
cluster tilting subcategory. Recall that the cluster tilting subcategories of $\cC(Z)$ have the form $\cT( T )$ where $T$ is a triangulation of $Z$.
\subsection{Existence of inf and sup}
\label{inf-sup}
\begin{figure}
  \centering
    \begin{tikzpicture}[scale=2.5]
      \draw (0,0) circle (1cm);
      
      \draw[very thick] ([shift=(-15:1cm)]0,0) arc (-15:55:1cm);      
      \draw[very thick] ([shift=(145:1cm)]0,0) arc (145:225:1cm);            

      \draw (-15:0.97cm) -- (-15:1.03cm);
      \draw (-15:1.13cm) node{$x_0$};
      \draw (8:0.97cm) -- (8:1.03cm);
      \draw (5.5:1.45cm) node{$\inf_{ \llbracket x_0,y_0 \rrbracket }M$};
      \draw (37:0.97cm) -- (37:1.03cm);
      \draw (27:1.43cm) node{$\sup_{ \llbracket x_0,y_0 \rrbracket }M$};
      \draw (55:0.97cm) -- (55:1.03cm);
      \draw (55:1.13cm) node{$y_0$};
      \draw (145:0.97cm) -- (145:1.03cm);
      \draw (145:1.13cm) node{$x_1$};
      \draw (170:0.97cm) -- (170:1.03cm);
      \draw (210:0.97cm) -- (210:1.03cm);
      \draw (225:0.97cm) -- (225:1.03cm);
      \draw (225:1.15cm) node{$y_1$};


      \draw[red] (8:1cm) .. controls (8:0.45cm) and (210:0.45cm) .. (210:1cm);
                  
      \draw[red] (37:1cm) .. controls (37:0.45cm) and (170:0.45cm) .. (170:1cm);

      \draw[red] (185:0.8cm) node{$\cdot$};
      \draw[red] (192.5:0.8cm) node{$\cdot$};
      \draw[red] (200:0.8cm) node{$\cdot$};

      \draw[red] (15:0.8cm) node{$\cdot$};
      \draw[red] (22.5:0.8cm) node{$\cdot$};
      \draw[red] (30:0.8cm) node{$\cdot$};

    \end{tikzpicture} 
  \caption{The set $M$ in Lemma \ref{lU} consists of the vertices in $\llbracket x_0,y_0 \rrbracket$ which are connected to $\llbracket x_1,y_1 \rrbracket$ by an arc in $T$.  The lemma shows that inf and sup of $M$ are vertices in $Z$.}
\label{fig:M}
\end{figure}
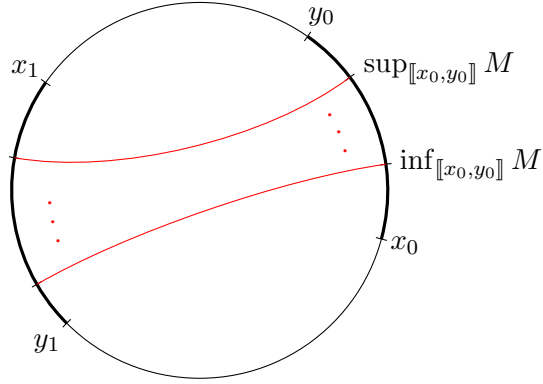

\ble{U}
Let $T$ be a triangulation of $Z$ and let $x_0 \leqslant y_0 < x_1 \leqslant y_1 < x_0$ be in $Z$.  Suppose the following set is non-empty:
\[
  M = \big\{ v \in Z \cap \llbracket x_0,y_0 \rrbracket \,\big|\, \mbox{$v$ is
    connected to $\llbracket x_1,y_1 \rrbracket$ by an edge or a diagonal in $T$} \big\},
\]
see Figure \ref{fig:M}.
Then $\inf_{ \llbracket x_0,y_0 \rrbracket }M$ and $\sup_{ \llbracket x_0,y_0 \rrbracket }M$ are in $Z$ and hence in $M$.
\ele

\begin{proof}
Suppose that $i = \inf_{ \llbracket x_0,y_0 \rrbracket }M$ is not in $Z$.  Then $i$ is a limit point of $Z$, and there are diagonals $\{ v_j,w_j \}$ in $T$ for $j \in \BN$ such that $v_j \in \llbracket x_0,y_0 \rrbracket$, $w_j \in \llbracket x_1,y_1 \rrbracket$, and $v_j \searrow i$.  We can assume $x_0 < i < \cdots < v_2 < v_1 < y_0$ whence the non-crossing condition on $T$ implies $x_1 \leqslant w_1 \leqslant w_2 \leqslant \cdots \leqslant y_1$.  We will show that there is neither a leapfrog nor a fountain in $T$ converging to $i$; this is a contradiction.

There is no leapfrog in $T$ converging to $i$:  Suppose there is such a leapfrog $\{ p_k,q_k \}$ for $k \in \BN$.  Pick $K$ such that $x_0 \leqslant p_K < i < q_K \leqslant y_0$.  Then pick $J$ such that $i < v_J < q_K$.  This implies $x_0 \leqslant p_K < i < v_J < q_K \leqslant y_0 < x_1 \leqslant w_J \leqslant y_1$, so in particular $p_K < v_J < q_K < w_J$, contradicting that $T$ is non-crossing.

There is no fountain in $T$ converging to $i$:  Suppose there is such a fountain $\{ p,q_k \}$ for $k \in \BN$.  Pick $K$ and $L$ such that $x_0 \leqslant q_K < i < q_L \leqslant y_0$.  Pick $J$ such that $i < v_J < q_L$.    The diagonal $\{ v_J,w_J \}$ does not cross $\{ p,q_K \}$ or $\{ p,q_L \}$, because these diagonals are all in $T$.  But we know $v_J \in \llbracket q_K,q_L \rrbracket$, so must have either $w_J = p$ or $w_J \in \llbracket q_K,q_L \rrbracket$.  The latter would imply $w_J \in \llbracket x_0,y_0 \rrbracket$ which is false, so the former must hold: $w_J = p$.  Then $p \in \llbracket x_1,y_1 \rrbracket$ whence $\{ p,q_K \} \in M$, but this contradicts the definition of $i$ because $x_0 \leqslant q_K < i$. 

The proof that $\sup_{ \llbracket x_0,y_0 \rrbracket }M$ is in $Z$ is symmetric.
\end{proof}
\subsection{The zig-zag path}
\label{zig-zag}

\begin{figure}
  \centering
    \begin{tikzpicture}[scale=4]
      \draw (0,0) circle (1cm);

      \draw (50:0.97cm) -- (50:1.03cm);
      \draw (50:1.13cm) node{$e_{ 2i-3 }$};
      \draw (75:0.97cm) -- (75:1.03cm);
      \draw (75:1.13cm) node{$e_{ 2i-1 }$};
      \draw (100:0.97cm) -- (100:1.03cm);
      \draw (100:1.13cm) node{$e_{ 2i+1 } = f$};
      \draw (125:0.97cm) -- (125:1.03cm);
      \draw (125:1.13cm) node{$e_{ 2i }$};
      \draw (150:0.97cm) -- (150:1.03cm);
      \draw (150:1.15cm) node{$e_{ 2i-2 }$};

      \draw (200:0.97cm) -- (200:1.03cm);
      \draw (200:1.13cm) node{$e_6$};
      \draw (220:0.97cm) -- (220:1.03cm);
      \draw (220:1.13cm) node{$e_4$};
      \draw (240:0.97cm) -- (240:1.03cm);
      \draw (240:1.13cm) node{$e_2$};
      \draw (260:0.97cm) -- (260:1.03cm);
      \draw (260:1.13cm) node{$e_0 = e$};
      \draw (280:0.97cm) -- (280:1.03cm);
      \draw (280:1.13cm) node{$e_1$};
      \draw (300:0.97cm) -- (300:1.03cm);
      \draw (300:1.13cm) node{$e_3$};
      \draw (320:0.97cm) -- (320:1.03cm);
      \draw (320:1.13cm) node{$e_5$};

      \draw (260:1cm) .. controls (265:0.85cm) and (275:0.85cm) .. (280:1cm);
      \draw (240:1cm) .. controls (250:0.75cm) and (270:0.75cm) .. (280:1cm);
      \draw (240:1cm) .. controls (250:0.65cm) and (290:0.65cm) .. (300:1cm);
      \draw (220:1cm) .. controls (230:0.55cm) and (290:0.55cm) .. (300:1cm);
      \draw (220:1cm) .. controls (230:0.45cm) and (310:0.45cm) .. (320:1cm);
      \draw (200:1cm) .. controls (210:0.35cm) and (310:0.35cm) .. (320:1cm);

      \draw (100:1cm) .. controls (105:0.85cm) and (120:0.85cm) .. (125:1cm);
      \draw (75:1cm) .. controls (85:0.75cm) and (115:0.75cm) .. (125:1cm);
      \draw (75:1cm) .. controls (85:0.65cm) and (140:0.65cm) .. (150:1cm);
      \draw (50:1cm) .. controls (60:0.55cm) and (140:0.55cm) .. (150:1cm);

      \draw (-0.05,-0.05) node{$\cdot$};
      \draw (-0.05,0.05) node{$\cdot$};
      \draw (-0.05,0.15) node{$\cdot$};

    \end{tikzpicture} 
  \caption{The zig-zag of \prref{pro:zig-zag}.}
\label{fig:zig-zag}
\end{figure}
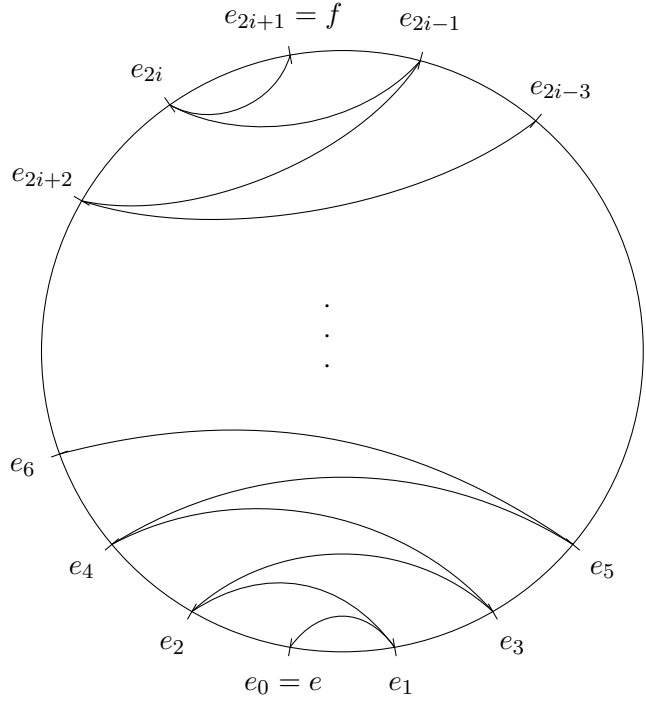

\bpr{pro:zig-zag}
Let $T$ be a triangulation of $Z$ and $e \neq f$ be in $Z$.  There exist an integer $i \geqslant 0$ and
a sequence of vertices $e_0, \ldots, e_{ 2i+1 }$ in $Z$ satisfying the
following, see Figure \ref{fig:zig-zag}.
\begin{enumerate}
\setlength\itemsep{4pt}

  \item  $e_0 = e$ and $e_{ 2i+1 } = f$.

  \item  If $\{ e_m, e_{ m+1 } \}$ is a diagonal, then it is in $T$.

  \item  $e_0 < e_1 < e_3 < e_5 < \cdots < f < \cdots < e_6 < e_4 < e_2 <e$.

  \item  The vertices $e_{ 2j+1 }$ are defined as follows:

\begin{enumerate}
\setlength\itemsep{4pt}

  \item  $e_1$ is the last vertex in $\llbracket e_0^+,f \rrbracket$ which is
          connected to $e_0$ by an edge or a diagonal in $T$.

  \item  If $j \geqslant 1$ then $e_{ 2j+1 }$ is the last vertex in $\llbracket e_{ 2j-1 }^+,f \rrbracket$ which is connected to $e_{ 2j }$ by an edge or a diagonal in $T$. 

\end{enumerate}

  \item  The vertices $e_{ 2j }$ are defined as follows:  If $j \geqslant 1$ then $e_{ 2j }$ is the first vertex in $\llbracket f^+,e_{2j-2 }^- \rrbracket$ which is connected to $e_{ 2j-1 }$ by a diagonal in $T$. 

\end{enumerate}
\epr

\begin{proof}
Set $e_0 = e$.

Observe that the set 
\[
  U = \big\{ u \in Z \cap \llbracket e_0^+,f \rrbracket \,\big|\, \mbox{$u$ is connected to $e_0$ by an edge or a diagonal in $T$} \big\} 
\]
is non-empty because $e_0^+ \in U$.  Set
\[
  e_1 = \sup_{\llbracket e_0^+,f \rrbracket} U.
\]
Applying the special case of \leref{U} when two of the four points coincide gives that $e_1$ is a vertex of $U$, which clearly satisfies (4)(a).  If $e_1 = f$ then terminate the construction.  Otherwise we have $e_0 < e_1 < f$; in particular, the sequence consisting of $e_0$ and $e_1$ satisfies (2)--(5).  

Let $k \geqslant 1$ be given, and suppose that the sequence of
vertices $e_0, \ldots, e_{ 2k-1 }$ has been defined such that
(2)--(5) are satisfied and such that $e_{ 2k-1 } \neq f$.  By (3)
we have $e_0 < e_{ 2k-1 } < f < e_{ 2k-2 }$.

Consider the set
\[
  V_{ 2k } = \big\{ v \in Z \cap \llbracket f^+,e_{ 2k-2 }^- \rrbracket \,\big|\, \mbox{$v$ is connected to $e_{ 2k-1 }$ by a diagonal in $T$} \big\}, 
\]
see Figure \ref{fig:V2k}.
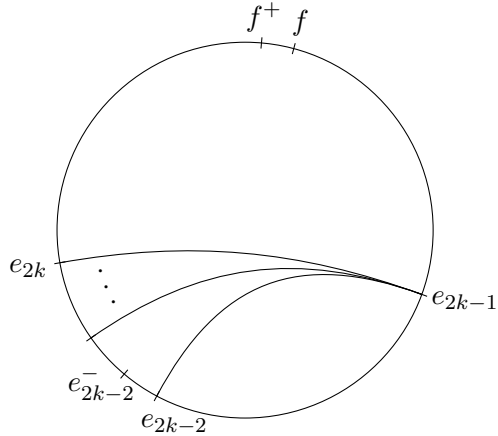
\begin{figure}
  \centering
    \begin{tikzpicture}[scale=2.5]
      \draw (0,0) circle (1cm);

      \draw (-20:0.97cm) -- (-20:1.03cm);
      \draw (-18:1.23cm) node{$e_{ 2k-1 }$};
      \draw (75:0.97cm) -- (75:1.03cm);
      \draw (75:1.14cm) node{$f$};
      \draw (85:0.97cm) -- (85:1.03cm);
      \draw (85:1.14cm) node{$f^+$};
      \draw (190:0.97cm) -- (190:1.03cm);
      \draw (190:1.17cm) node{$e_{ 2k }$};
      \draw (215:0.97cm) -- (215:1.03cm);
      \draw (230:0.97cm) -- (230:1.03cm);
      \draw (227:1.12cm) node{$e_{ 2k-2 }^-$};
      \draw (242:0.97cm) -- (242:1.03cm);
      \draw (250:1.10cm) node{$e_{ 2k-2 }$};

      \draw (-20:1cm) .. controls (-20:0.3cm) and (190:0.3cm) .. (190:1cm);
      \draw (-20:1cm) .. controls (-20:0.3cm) and (215:0.3cm) .. (215:1cm);
      \draw (-20:1cm) .. controls (-20:0.3cm) and (242:0.3cm) .. (242:1cm);

      \draw (196:0.8cm) node{$\cdot$};
      \draw (202.5:0.8cm) node{$\cdot$};
      \draw (209:0.8cm) node{$\cdot$};

    \end{tikzpicture} 
  \caption{The set $V_{ 2k }$ consists of the vertices in 
  $\llbracket f^+,e_{2k-2 }^- \rrbracket$ which are connected to $e_{ 2k-1 }$ by a diagonal in $T$.} 
\label{fig:V2k}
\end{figure}
Observe that $V_{ 2k }$ is non-empty: If it was empty, then since $T$
is a triangulation, $e_{ 2k-2 }$ would be connected by a diagonal in $T$
to a vertex $z \in \llbracket e_{ 2k-1 }^+,f \rrbracket$, contradicting (4)(b).
Set
\[
  e_{ 2k } = \inf_{ \llbracket f^+,e_{ 2k-2 }^-  \rrbracket} V_{ 2k }.
\]
Again applying the special case of \leref{U} gives that $e_{ 2k }$ is a vertex of $V_{ 2k }$, which clearly satisfies (5).  We have
\begin{equation}
\label{equ:inequality_even}
  f < e_{ 2k } < e_{ 2k-2 }.
\end{equation}

Now consider the set
\[
  W_{ 2k+1 } = \big\{ w \in Z \cap \llbracket e_{ 2k-1 }^+,f \rrbracket \,\big|\, \mbox{$w$ is connected to $e_{ 2k }$ by an edge or a diagonal in $T$} \big\},
\]
see Figure \ref{fig:W2k1}.
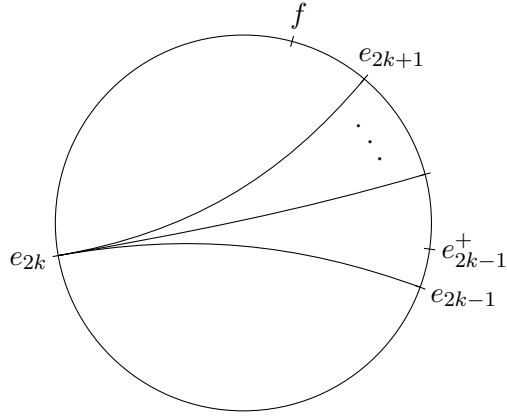
\begin{figure}
  \centering
    \begin{tikzpicture}[scale=2.5]
      \draw (0,0) circle (1cm);

      \draw (190:0.97cm) -- (190:1.03cm);
      \draw (190:1.16cm) node{$e_{ 2k }$};
      \draw (-20:0.97cm) -- (-20:1.03cm);
      \draw (-19:1.24cm) node{$e_{ 2k-1 }$};
      \draw (-8:0.97cm) -- (-8:1.03cm);
      \draw (-7:1.24cm) node{$e_{ 2k-1 }^+$};
      \draw (15:0.97cm) -- (15:1.03cm);
      \draw (50:0.97cm) -- (50:1.03cm);
      \draw (48:1.18cm) node{$e_{ 2k+1 }$};
      \draw (75:0.97cm) -- (75:1.03cm);
      \draw (75:1.14cm) node{$f$};

      \draw (190:1cm) .. controls (190:0.3cm) and (-20:0.3cm) .. (-20:1cm);
      \draw (190:1cm) .. controls (190:0.3cm) and (12:0.3cm) .. (15:1cm);
      \draw (190:1cm) .. controls (190:0.3cm) and (50:0.3cm) .. (50:1cm);

      \draw (25:0.8cm) node{$\cdot$};
      \draw (32.5:0.8cm) node{$\cdot$};
      \draw (40:0.8cm) node{$\cdot$};

    \end{tikzpicture} 
  \caption{The set $W_{ 2k+1 }$ consists of the vertices in 
      $\llbracket e_{2k-1 }^+,f \rrbracket$ which are connected to $e_{ 2k }$ by an edge or a diagonal in $T$.} 
\label{fig:W2k1}
\end{figure}
Observe that $W_{ 2k+1 }$ is non-empty: If it was empty, then since $T$
is a triangulation, $e_{ 2k-1 }$ would be connected by a diagonal in $T$
to a vertex $z \in \llbracket f^+,e_{ 2k }^- \rrbracket$, contradicting (5).  Set
\[
  e_{ 2k+1 } = \sup_{ \llbracket e_{ 2k-1 }^+,f \rrbracket } W_{ 2k+1 }.
\]
\leref{U} says that $e_{ 2k+1 }$ is a vertex of $W_{ 2k+1 }$ which
clearly satisfies (4)(b).  If $e_{ 2k+1 } = f$ then terminate the
construction.  Otherwise we have
\begin{equation}
\label{equ:inequality_odd}
  e_{ 2k-1 } < e_{ 2k+1 } < f.
\end{equation}

We have now extended the sequence of vertices to $e_0, \ldots, e_{
  2k+1 }$ and shown that the extended sequence satisfies (4) and (5).
It is clear that it satisfies (2), and (3) holds by Equations
\eqref{equ:inequality_even} and \eqref{equ:inequality_odd}.

To complete the proof, we must show that the construction terminates.
If not, then (3) would mean that $e_{ 2j +1 } \nearrow p$ and $e_{
  2j } \searrow q$ for certain $p$, $q$ in $\overline{Z}$ with $p
\leqslant f \leqslant q$.

If we had $p = q$ then we would get $p = f = q$, but this would mean
$e_{ 2j+1 } = f$ for some $j$, contradicting that the construction
did not terminate.

If we had $p \neq q$, then note that \cite[thm.\ 0.5 and prop.\ 5.6]{GHJ} implies
that $T$ satisfies condition (PC2) of \cite{GHJ}.  Hence
there would be sequences $\{ x_{ \ell } \}$ and $\{ y_{ \ell } \}$ in
$Z$ such that $\{ x_{ \ell },y_{ \ell } \} \in T$ for each $\ell$, and
such that $x_{ \ell } \searrow p$ and $y_{ \ell } \searrow q$.  But
then it is easy to see that for some $j$ and some $\ell$, the diagonals
$\{ e_{ 2j },e_{ 2j+1 } \}$ and $\{ x_{ \ell },y_{ \ell } \}$ would
cross, contradicting that they are both in the triangulation $T$.
\end{proof}
\subsection{The index formula}
\label{ind-cC}
We can now prove the formula for the index of an indecomposable object of $\cC( Z )$ with respect to an arbitrary cluster tilting subcategory.
\bth{ind-formula}
Let $e \neq f$ be in $Z$ and and $T$ be a triangulation of $Z$. Consider the vertices $e = e_0, \ldots,
e_{ 2i+1 } = f$ of the zig-zag path from \prref{pro:zig-zag}, see Figure
\ref{fig:zig-zag}.  Then
\[
  \ind_{ \cT } \{ e,f \}
  = \sum_{ m = 0 }^{ 2i } (-1)^m \big[ \{ e_m,e_{ m+1 } \} \big].
\]
\eth

\begin{proof}
The theorem is trivial for $i = 0$.  Let $i \geqslant 1$ be given and
suppose the theorem holds for lower values of $i$.  \prref{pro:zig-zag} gives the crossing diagonals in Figure \ref{fig:crossing},
\begin{figure}
  \centering
    \begin{tikzpicture}[scale=2.5]
      \draw (0,0) circle (1cm);

      \draw (280:0.97cm) -- (280:1.03cm);
      \draw (280:1.17cm) node{$e_0$};
      \draw (70:0.97cm) -- (70:1.03cm);
      \draw (70:1.20cm) node{$e_{ 2i-1 }$};
      \draw (120:0.97cm) -- (120:1.03cm);
      \draw (120:1.20cm) node{$e_{ 2i+1 }$};
      \draw (160:0.97cm) -- (160:1.03cm);
      \draw (160:1.17cm) node{$e_{ 2i }$};

      \draw (70:1cm) .. controls (95:0.5cm) and (135:0.5cm) .. (160:1cm);
      \draw[red] (280:1cm) .. controls (240:0.1cm) and (160:0.1cm) .. (120:1cm);
      \draw[blue] (280:1cm) .. controls (290:0.1cm) and (60:0.1cm) .. (70:1cm);
      \draw[blue] (120:1cm) .. controls (130:0.7cm) and (150:0.7cm) .. (160:1cm);

    \end{tikzpicture} 
  \caption{A crossing resulting from \prref{pro:zig-zag}.}
\label{fig:crossing}
\end{figure}
hence the following triangle in $\cC( Z )$ by \cite[p.\ 4389]{IT2}.
\begin{equation}
\label{equ:crossing}
  \{ e_{ 2i-1 },e_{ 2i } \}
  \rightarrow
  \{ e_{ 2i },e_{ 2i+1 } \} \oplus \{ e_0,e_{ 2i-1 } \}
  \stackrel{ g }{ \rightarrow }
  \{ e_0,e_{ 2i+1 } \}
  \rightarrow
  \Sigma \{ e_{ 2i-1 },e_{ 2i } \}
\end{equation}
The functor
\[
  F : \cC( Z ) \rightarrow \mod\,\cT
  \;\;,\;\;
  Fc = \Hom_{ \cC( Z ) }( -,c )|_{ \cT }
\]
sends \eqref{equ:crossing} to an exact sequence, and it sends $\Sigma
\{ e_{ 2i-1 },e_{ 2i } \}$ to $0$ because $\{ e_{ 2i-1 },e_{ 2i } \}
\in \cT$.  Hence $Fg$ is surjective, so \cite[lem.\ 2.2]{P} says that
$\ind_{ \cT }$ is additive on the extension given by
\eqref{equ:crossing}, that is,
\[
  \ind_{ \cT } \big( \{ e_{ 2i },e_{ 2i+1 } \} \oplus \{ e_0,e_{ 2i-1 } \} \big)
  =
  \ind_{ \cT } \{ e_{ 2i-1 },e_{ 2i } \}
  + 
  \ind_{ \cT } \{ e_0,e_{ 2i+1 } \}.
\]
Rearranging gives
\[
  \ind_{ \cT } \{ e_0,e_{ 2i+1 } \}
  =
  \ind_{ \cT } \{ e_0,e_{ 2i-1 } \}
  + \big[ \{ e_{ 2i },e_{ 2i+1 } \} \big]
  - \big[ \{ e_{ 2i-1 },e_{ 2i } \} \big],
\]
and combining with the theorem for $i-1$ proves the theorem for $i$.
\end{proof}
\sectionnew{Classification of the $c$-vectors of the categories $\cC(Z)$}
\lb{icvect}
In this section we give an explicit classification of the set of positive $c$-vectors of the categories $\cC(Z)$ with 
respect to an arbitrary cluster tilting subcategory $\cT(T)$.  The result appears in \thref{c-vectors_and_virtual_dimension_vectors}. Along the way we 
describe in an explicit way the images in $\mod\big( \cT( T ) \big)$ 
that appear in \thref{2}(2). We furthermore show that all such images are indecomposable objects of $\mod\big( \cT( T ) \big)$, and, in the 
opposite direction, that the dimension vectors of all such indecomposables are positive $c$-vectors of $\cC(Z)$. 
\subsection{The triangles of $T$} 
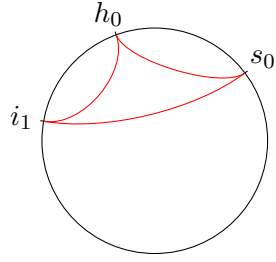
\begin{figure}
  \centering
    \begin{tikzpicture}[scale=1.5]
      \draw (0,0) circle (1cm);

      \draw (37:0.97cm) -- (37:1.03cm);
      \draw (37:1.20cm) node{$s_0$};
      \draw (110:0.97cm) -- (110:1.03cm);
      \draw (110:1.20cm) node{$h_0$};      
      \draw (170:0.97cm) -- (170:1.03cm);
      \draw (170:1.20cm) node{$i_1$};
                  
      \draw[red] (37:1cm) .. controls (37:0.45cm) and (170:0.45cm) .. (170:1cm);
      \draw[red] (37:1cm) .. controls (37:0.75cm) and (110:0.75cm) .. (110:1cm);
      \draw[red] (110:1cm) .. controls (110:0.65cm) and (170:0.65cm) .. (170:1cm);

    \end{tikzpicture} 
  \caption{Let $\{ i_1,s_0 \}$ be a diagonal in $T$.  Lemma \ref{ltriangles} shows that there is a triangle of $T$ as shown.}
\label{fig:triangles}
\end{figure}

We start with three lemmas on the combinatorics of triangles adjacent to diagonals in $T$. 
\ble{triangles}
Let $T$ be a triangulation of $Z$ and suppose that $\{ s_0,i_1 \}$ is a diagonal in $T$.  Then there exists $h_0 \in Z$ such that $s_0 < h_0 < i_1$ and such that $\{ i_1,h_0,s_0 \}$ is a triangle of $T$, see Figure \ref{fig:triangles}.

\ele

\begin{proof}
Apply Lemma \ref{lU} with $x_0 = s_0^+$, $y_0 = i_1^-$, $x_1 = y_1 = i_1$.  The set $M$ is non-empty because it contains the edge $\{ i_1^-,i_1 \}$, and $h_0 = \inf_{ \llbracket s_0^+,i_1^- \rrbracket }M$ is in $Z$ by Lemma \ref{lU}.  Then it is clear that $s_0 < h_0 < i_1$ and that $\{ h_0,i_1 \}$ is an edge or diagonal in $T$.  It is easy to check that $T$ cannot contain a diagonal crossing $\{ h_0,s_0 \}$ whence $\{ h_0,s_0 \}$ is an edge or a diagonal in $T$. 
\end{proof}

\begin{figure}
  \centering
    \begin{tikzpicture}[scale=2.5]
      \draw (0,0) circle (1cm);

      \draw[very thick] ([shift=(-15:1cm)]0,0) arc (-15:55:1cm);      
      \draw[very thick] ([shift=(145:1cm)]0,0) arc (145:225:1cm);            

      \draw (-15:0.97cm) -- (-15:1.03cm);
      \draw (-15:1.13cm) node{$a_0$};
      \draw (8:0.97cm) -- (8:1.03cm);
      \draw (8:1.13cm) node{$i_0$};
      \draw (37:0.97cm) -- (37:1.03cm);
      \draw (37:1.13cm) node{$s_0$};
      \draw (55:0.97cm) -- (55:1.03cm);
      \draw (55:1.13cm) node{$b_0$};
      \draw (110:0.97cm) -- (110:1.03cm);
      \draw (110:1.13cm) node{$h_0$};      
      \draw (145:0.97cm) -- (145:1.03cm);
      \draw (145:1.13cm) node{$a_1$};
      \draw (170:0.97cm) -- (170:1.03cm);
      \draw (170:1.13cm) node{$i_1$};
      \draw (210:0.97cm) -- (210:1.03cm);
      \draw (210:1.13cm) node{$s_1$};
      \draw (225:0.97cm) -- (225:1.03cm);
      \draw (225:1.15cm) node{$b_1$};
      \draw (265:0.97cm) -- (265:1.03cm);
      \draw (265:1.15cm) node{$h_1$};

      \draw[red] (8:1cm) .. controls (8:0.45cm) and (210:0.45cm) .. (210:1cm);
      \draw[red] (8:1cm) .. controls (8:0.60cm) and (265:0.60cm) .. (265:1cm);
      \draw[red] (210:1cm) .. controls (210:0.70cm) and (265:0.70cm) .. (265:1cm);
                  
      \draw[red] (37:1cm) .. controls (37:0.45cm) and (170:0.45cm) .. (170:1cm);
      \draw[red] (37:1cm) .. controls (37:0.75cm) and (110:0.75cm) .. (110:1cm);
      \draw[red] (110:1cm) .. controls (110:0.65cm) and (170:0.65cm) .. (170:1cm);

%

    \end{tikzpicture} 
  \caption{Let $a_0^- < b_0 < a_1^- < b_1 < a_0^-$ be in $Z$.  Suppose $\llbracket a_0,b_0 \rrbracket$ is connected to $\llbracket a_1,b_1 \rrbracket$ by at least one diagonal in $T$.  Then there exist edges or diagonals in $T$ as shown in red.  See Lemma \ref{lU_improved}.}
\label{fig:U_improved}
\end{figure}
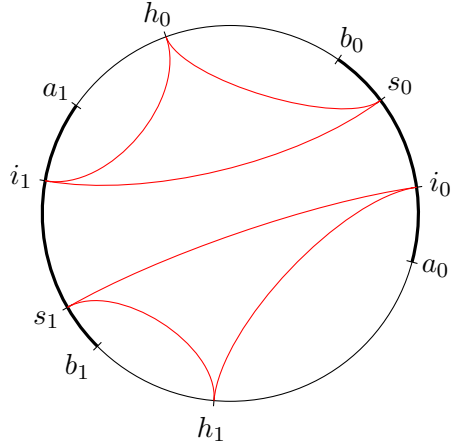

\ble{U_improved}
Let $T$ be a triangulation of $Z$, let $a_0^- < b_0 < a_1^- < b_1 < a_0^-$ be in $Z$, and suppose $\llbracket a_0,b_0 \rrbracket$ is connected to $\llbracket a_1,b_1 \rrbracket$ by at least one diagonal in $T$.  There exist
\[
  a_0 \leqslant i_0 \leqslant s_0 \leqslant b_0 < h_0 < a_1 \leqslant i_1 \leqslant s_1 \leqslant b_1 < h_1 < a_0
\]  
in $Z$ such that $\{ i_0,h_1,s_1 \}$ and $\{ i_1,h_0,s_0 \}$ are triangles of  $T$, see Figure \ref{fig:U_improved}.
\ele

\begin{proof}
Apply Lemma \ref{lU} with $x_0,y_0,x_1,y_1$ equal to $a_0,b_0,a_1,b_1$ and set $i_0 = \inf_{ \llbracket a_0,b_0 \rrbracket }M$.  Then apply Lemma \ref{lU} with $x_0,y_0,x_1,y_1$ equal to $a_1,b_1,i_0,i_0$ and set $s_1 = \sup_{ \llbracket a_1,b_1 \rrbracket }M$.  Symmetrically, apply Lemma \ref{lU} with $x_0,y_0,x_1,y_1$ equal to $a_1,b_1,a_0,b_0$ and set $i_1 = \inf_{ \llbracket a_1,b_1 \rrbracket }M$.  Then apply Lemma \ref{lU} with $x_0,y_0,x_1,y_1$ equal to $a_0,b_0,i_1,i_1$ and set $s_0 = \sup_{ \llbracket a_0,b_0 \rrbracket }M$.  It is easy to check that
\[
  a_0 \leqslant i_0 \leqslant s_0 \leqslant b_0 < a_1 \leqslant i_1 \leqslant s_1 \leqslant b_1.
\]
It also follows that $\{ i_1,s_0 \} \in T$, so by Lemma \ref{ltriangles} there is a triangle $\{ i_1,h_0,s_0 \}$ of $Z$ with $s_0 < h_0 < i_1$.  In fact we must have
\[
  b_0 < h_0 < a_1.
\]  
Namely, if $s_0 < h_0 \leqslant b_0$ then $\{ i_1,h_0 \} \in T$ would contradict the definition of $s_0$ as a supremum, and if $a_1 \leqslant h_0 < i_1$ then $\{ h_0,s_0 \} \in T$ would contradict the definition of $i_1$ as an infimum.  Symmetrically there is a triangle $\{ i_0,h_1,s_1 \}$ of $Z$ with
\[
  b_1 < h_1 < a_0.
\]
\end{proof}

We record the following special case of Lemma \ref{lU_improved} for later use.

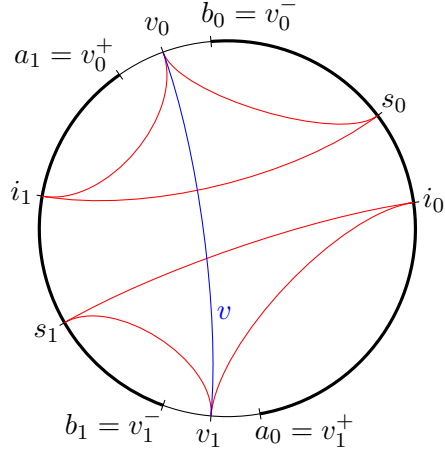
\begin{figure}
  \centering
    \begin{tikzpicture}[scale=2.5]
      \draw (0,0) circle (1cm);

      \draw[very thick] ([shift=(-80:1cm)]0,0) arc (-80:95:1cm);      
      \draw[very thick] ([shift=(125:1cm)]0,0) arc (125:250:1cm);            

      \draw (-80:0.97cm) -- (-80:1.03cm);
      \draw (-69:1.15cm) node{$a_0 = v_1^+$};
      \draw (8:0.97cm) -- (8:1.03cm);
      \draw (8:1.11cm) node{$i_0$};
      \draw (37:0.97cm) -- (37:1.03cm);
      \draw (37:1.11cm) node{$s_0$};
      \draw (95:0.97cm) -- (95:1.03cm);
      \draw (84:1.15cm) node{$b_0 = v_0^-$};
      \draw (110:0.97cm) -- (110:1.03cm);
      \draw (109:1.14cm) node{$v_0$};      
      \draw (125:0.97cm) -- (125:1.03cm);
      \draw (133:1.27cm) node{$a_1 = v_0^+$};
      \draw (170:0.97cm) -- (170:1.03cm);
      \draw (170:1.11cm) node{$i_1$};
      \draw (210:0.97cm) -- (210:1.03cm);
      \draw (210:1.11cm) node{$s_1$};
      \draw (250:0.97cm) -- (250:1.03cm);
      \draw (240:1.21cm) node{$b_1 = v_1^-$};
      \draw (265:0.97cm) -- (265:1.03cm);
      \draw (265:1.12cm) node{$v_1$};

      \draw[red] (8:1cm) .. controls (8:0.45cm) and (210:0.45cm) .. (210:1cm);
      \draw[red] (8:1cm) .. controls (8:0.60cm) and (265:0.60cm) .. (265:1cm);
      \draw[red] (210:1cm) .. controls (210:0.70cm) and (265:0.70cm) .. (265:1cm);
                  
      \draw[red] (37:1cm) .. controls (37:0.45cm) and (170:0.45cm) .. (170:1cm);
      \draw[red] (37:1cm) .. controls (37:0.75cm) and (110:0.75cm) .. (110:1cm);
      \draw[red] (110:1cm) .. controls (110:0.65cm) and (170:0.65cm) .. (170:1cm);
      
      \draw[blue] (110:1cm) .. controls (110:0.45cm) and (265:0.45cm) .. (265:1cm);

      \draw[blue] (268:0.45cm) node{$v$};

    \end{tikzpicture} 
  \caption{Let $v$ be a diagonal of $Z$ crossing at least one diagonal in $T$.  Then there exist edges or diagonals in $T$ as shown in red.  See Lemma \ref{lU_improved_special_case}.}
\label{fig:U_improved_special_case}
\end{figure}

\ble{U_improved_special_case}
Let $T$ be a triangulation of $Z$ and let $v = \{ v_0,v_1 \}$ be a diagonal of $Z$ which crosses at least one diagonal in $T$.  There exist
\[
  i_0 \leqslant s_0 < v_0 < i_1 \leqslant s_1 < v_1
\]
in $Z$  such that $\{ i_0,v_1,s_1 \}$ and $\{ i_1,v_0,s_0 \}$ are triangles of $T$, see Figure \ref{fig:U_improved_special_case}.
\ele

\begin{proof}
Set $a_0 = v_1^+$, $b_0 = v_0^-$, $a_1 = v_0^+$, $b_1 = v_1^-$.  The diagonals in $T$ which cross $v$ are precisely the diagonals in $T$ which connect $\llbracket a_0,b_0 \rrbracket$ to $\llbracket a_1,b_1 \rrbracket$.  We can apply Lemma \ref{lU_improved}, and this proves the present lemma because we must have $h_0 = v_0$, $h_1 = v_1$.
\end{proof}

\subsection{A formula for images in $\mod\big( \cT( T ) \big)$} Our next result provides an explicit description of the images in $\mod\big( \cT( T ) \big)$ which appear in \thref{2}(2). 

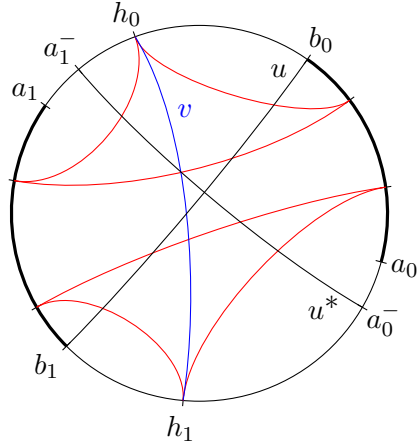
\begin{figure}
  \centering
    \begin{tikzpicture}[scale=2.5]
      \draw (0,0) circle (1cm);

      \draw[very thick] ([shift=(-15:1cm)]0,0) arc (-15:55:1cm);      
      \draw[very thick] ([shift=(145:1cm)]0,0) arc (145:225:1cm);            

      \draw (-30:0.97cm) -- (-30:1.03cm);
      \draw (-30:1.13cm) node{$a^-_0$};      
      \draw (-15:0.97cm) -- (-15:1.03cm);
      \draw (-15:1.13cm) node{$a_0$};
      \draw (8:0.97cm) -- (8:1.03cm);
      \draw (37:0.97cm) -- (37:1.03cm);
      \draw (55:0.97cm) -- (55:1.03cm);
      \draw (55:1.13cm) node{$b_0$};
      \draw (110:0.97cm) -- (110:1.03cm);
      \draw (110:1.13cm) node{$h_0$};      
      \draw (130:0.97cm) -- (130:1.03cm);
      \draw (130:1.15cm) node{$a_1^-$};      
      \draw (145:0.97cm) -- (145:1.03cm);
      \draw (145:1.13cm) node{$a_1$};
      \draw (170:0.97cm) -- (170:1.03cm);
      \draw (210:0.97cm) -- (210:1.03cm);
      \draw (225:0.97cm) -- (225:1.03cm);
      \draw (225:1.15cm) node{$b_1$};
      \draw (265:0.97cm) -- (265:1.03cm);
      \draw (265:1.15cm) node{$h_1$};

      \draw[red] (8:1cm) .. controls (8:0.45cm) and (210:0.45cm) .. (210:1cm);
      \draw[red] (8:1cm) .. controls (8:0.60cm) and (265:0.60cm) .. (265:1cm);
      \draw[red] (210:1cm) .. controls (210:0.70cm) and (265:0.70cm) .. (265:1cm);
                  
      \draw[red] (37:1cm) .. controls (37:0.45cm) and (170:0.45cm) .. (170:1cm);
      \draw[red] (37:1cm) .. controls (37:0.75cm) and (110:0.75cm) .. (110:1cm);
      \draw[red] (110:1cm) .. controls (110:0.65cm) and (170:0.65cm) .. (170:1cm);

      \draw[blue] (110:1cm) .. controls (95:0.45cm) and (270:0.45cm) .. (265:1cm);      
      \draw[blue] (98:0.55cm) node{$v$};
      
      \draw (-30:1cm) .. controls (-30:0.45cm) and (130:0.45cm) .. (130:1cm);      
      \draw (55:1cm) .. controls (55:0.45cm) and (225:0.45cm) .. (225:1cm);      
      \draw (61:0.87cm) node{$u$};
      \draw (-37:0.83cm) node{$u^{ \ast }$};            
      
%

    \end{tikzpicture} 
  \caption{There is a non-zero morphism $u \rightarrow \Sigma u^{ \ast }$ which induces a morphism
$
    \varphi : \Hom_{ \cC( Z ) }( -,u )\big|_{ \cT( T ) }
              \rightarrow
              \Hom_{ \cC( Z ) }( -,\Sigma u^{ \ast } )\big|_{ \cT( T ) }
$.
If $\Image \varphi \neq 0$ then $\Image \varphi = \Hom_{ \cC( Z ) }( -,\Sigma v )\big|_{ \cT( T ) }$, where $v = \{ h_0,h_1 \}$ is determined by the existence of red edges or diagonals in $T$ as shown.  See Proposition \ref{psubobjects}.}
\label{fig:image}
\end{figure}

\bpr{subobjects}
Let $T$ be a triangulation of $Z$. Let $u \rightarrow \Sigma u^{ \ast }$ be a non-zero morphism
with $u,u^{ \ast } \in \Indec\, \cC( Z )$ and consider the induced morphism
\begin{equation}
\label{equ:phi}
  \varphi : \Hom_{ \cC( Z ) }( -,u )\big|_{ \cT( T ) }
            \rightarrow
            \Hom_{ \cC( Z ) }( -,\Sigma u^{ \ast } )\big|_{ \cT( T ) }.
\end{equation}
The diagonals of $u$ and $u^{ \ast }$ cross, so we can pick
\begin{equation}
\label{equ:vertex_inequalities}
  b_0 < a_1^- < b_1 < a_0^- < b_0
\end{equation}
in $Z$ such that $u = \{ b_0,b_1 \}$, $u^{ \ast } = \{ a_0^-,a_1^- \}$, see Figure \ref{fig:image}.
\begin{enumerate}
\setlength\itemsep{4pt}

  \item  If $\llbracket a_0,b_0 \rrbracket$ is not connected to $\llbracket a_1,b_1 \rrbracket$ by a diagonal in $T$, then $\Image \varphi = 0$.

  \item  If $\llbracket a_0,b_0 \rrbracket$ is connected to $\llbracket a_1,b_1 \rrbracket$ by at least one diagonal in $T$, then Lemma \ref{lU_improved} provides $h_0,h_1 \in Z$, and
\begin{equation}
\label{equ:Image}
  \Image \varphi = \Hom_{ \cC( Z ) }( -,\Sigma v )\big|_{ \cT( T ) }
\end{equation} 
where $v = \{ h_0,h_1 \}$.  In particular, $\Image \varphi$ is indecomposable.

\end{enumerate}

\epr

\begin{proof}
We have $u = \{ b_0,b_1 \}$ and $\Sigma u^{ \ast } = \{ a_0^{ -- },a_1^{ -- } \}$.  Equation \eqref{equ:vertex_inequalities} implies $b_0 \leqslant a_1^{ -- } < b_1 \leqslant a_0^{ -- } < b_0$.  If $t \in T$ then by Remark \ref{rIgusa-Todorov_morphisms}(2) there is a non-zero composition $t \rightarrow u \rightarrow \Sigma u^{ \ast }$ if and only if $t = \{ t_0,t_1 \}$ with $t_0 \leqslant b_0 \leqslant a_1^{ -- } \leqslant t_1^{ -- } < t_1 \leqslant b_1 \leqslant a_0^{ -- } \leqslant t_0^{ -- }$.  In particular, this can only happen if $a_0 \leqslant t_0 \leqslant b_0$ and $a_1 \leqslant t_1 \leqslant b_1$, that is, if $\llbracket a_0,b_0 \rrbracket$ is connected to $\llbracket a_1,b_1 \rrbracket$ by at least one diagonal in $T$, so part (1) of the proposition follows.

For part (2), assume that $\llbracket a_0,b_0 \rrbracket$ is connected to $\llbracket a_1,b_1 \rrbracket$ by at least one diagonal in $T$.  Lemma \ref{lU_improved} provides the vertices $i_0,s_0,h_0,i_1,s_1,h_1 \in Z$ and we have
\[
  b_0 \leqslant h_0^- \leqslant a_1^{ -- } \leqslant b_1^{ -- } <
  b_1 \leqslant h_1^- \leqslant a_0^{ -- } \leqslant b_0^{ -- }.
\]  
Since
\[
  u = \{ b_0,b_1 \}
  \;\;,\;\;
  \Sigma v = \{ h_0^-,h_1^- \}
  \;\;,\;\;
  \Sigma u^{ \ast } = \{ a_0^{ -- },a_1^{ -- } \},
\]  
Remark \ref{rIgusa-Todorov_morphisms} says that there are morphisms
\[
  u \rightarrow \Sigma v \rightarrow \Sigma u^{ \ast }
\]  
with a non-zero composition.  Up to a scalar multiple, the composition is the non-zero morphism $u \rightarrow \Sigma u^{ \ast }$ from the proposition, so we can assume that the composition of the induced morphisms
\[
  \Hom_{ \cC( Z ) }( -,u )\big|_{ \cT( T ) }
  \stackrel{ \alpha }{ \longrightarrow }
  \Hom_{ \cC( Z ) }( -,\Sigma v )\big|_{ \cT( T ) }
  \stackrel{ \beta }{ \longrightarrow }
  \Hom_{ \cC( Z ) }( -,\Sigma u^{ \ast } )\big|_{ \cT( T ) }
\]
is $\varphi$.  To prove Equation \eqref{equ:Image} it is enough to show that $\alpha$ is an epimorphism, $\beta$ a monomorphism.  First some preparation:

Let $t \rightarrow \Sigma v$ be a non-zero morphism from $t = \{ t_0,t_1 \} \in T$ to $\Sigma v = \{ h_0^-,h_1^- \}$.  Remark \ref{rIgusa-Todorov_morphisms} says we can assume $t_0 \leqslant h_0^- \leqslant t_1^{ -- } < t_1 \leqslant h_1^- \leqslant t_0^{ -- }$ which implies $h_1^+ \leqslant t_0 \leqslant h_0^-$ and $h_0^+ \leqslant t_1 \leqslant h_1^-$.  This implies $i_0 \leqslant t_0 \leqslant s_0$ and $i_1 \leqslant t_1 \leqslant s_1$ because $\{ i_1,h_0,s_0 \}$ and $\{ i_0,h_1,s_1 \}$ are triangles of $T$ by Lemma \ref{lU_improved}, so in particular
\[
  a_0 \leqslant t_0 \leqslant b_0
  \;\; \mbox{ and } \;\;
  a_1 \leqslant t_1 \leqslant b_1,
\]
see Figure \ref{fig:U_improved}.  Combining with Lemma \ref{lU_improved} shows 
\begin{equation}
\label{equ:t_in_intervals}
  t_0
  \leqslant b_0
  \leqslant h_0^-
  \leqslant a_1^{ -- }
  \leqslant t_1^{ -- }
  < t_1
  \leqslant b_1
  \leqslant h_1^-
  \leqslant a_0^{ -- }
  \leqslant t_0^{ -- }.
\end{equation}

The morphism $\alpha$ is an epimorphism:  Equation \eqref{equ:t_in_intervals} implies
\[
  t_0 \leqslant b_0 \leqslant h_0^- \leqslant t_1^{ -- } < t_1 \leqslant b_1 \leqslant h_1^- \leqslant t_0^{ -- }.
\]
Since
\[
  t = \{ t_0,t_1 \}
  \;\;,\;\;
  u = \{ b_0,b_1 \}
  \;\;,\;\;
  \Sigma v = \{ h_0^-,h_1^- \},
\]
Remark \ref{rIgusa-Todorov_morphisms}(2) says that each morphism $t \rightarrow \Sigma v$ factorises as $t \rightarrow u \rightarrow \Sigma v$.  Hence $\alpha$ is an epimorphism. 

The morphism $\beta$ is a monomorphism:  Equation \eqref{equ:t_in_intervals} implies
\[
  t_0 \leqslant h_0^- \leqslant a_1^{ -- } \leqslant t_1^{ -- } < t_1 \leqslant h_1^- \leqslant a_0^{ -- } \leqslant t_0^{ -- }.
\]
Since
\[
  t = \{ t_0,t_1 \}
  \;\;,\;\;
  \Sigma v = \{ h_0^-,h_1^- \}
  \;\;,\;\;
  \Sigma u^{ \ast } = \{ a_0^{ -- },a_1^{ -- } \},
\]
Remark \ref{rIgusa-Todorov_morphisms}(2) says that the composition $t \rightarrow \Sigma v \rightarrow \Sigma u^{ \ast }$ is non-zero.  Hence $\beta$ is a monomorphism.
\end{proof}
\subsection{Positive $c$-vectors of $\cC(Z)$ are dimension vectors of indecomposables} 
We next use the results from the previous subsection to show that the positive  $c$-vectors of the categories $\cC(Z)$
are dimension vectors of indecomposable objects of $\mod( \cT(T) )$.
\bpr{c-vectors_are_dimension_vectors_of_indecomposables}
Let $\cT,\cU \subseteq \cC( Z )$ be cluster tilting subcategories and let $u \in \cU$ be indecomposable.

If $c_{ \cT }( u,\cU ) \geqslant 0$ then there exists an indecomposable object $m \in \mod( \cT )$ such that
\[
  c_{ \cT }( u,\cU ) = \dim_{ \cT }( m ). 
\]
\epr

\begin{proof}
Let $u^{ \ast }$ be the mutation of $u$ with respect to $\cU$ and let $u^{ \ast } \rightarrow e \rightarrow u \stackrel{ \delta }{ \rightarrow } \Sigma u^{ \ast }$ be an exchange triangle with respect to $\cU$.  The induced morphism 
\[
  \varphi :
    \Hom_{ \cC( Z ) }( -,u )\big|_{ \cT }
    \rightarrow
    \Hom_{ \cC( Z ) }( -,\Sigma u^{ \ast } )\big|_{ \cT }
\]
satisfies 
\[
  c_{ \cT }( u,\cU ) = \dim_{ \cT } ( \Image \varphi )
\]
by Theorem \ref{t2}(2), which applies to $\cC( Z )$ because the quiver of each cluster tilting subcategory of $\cC( Z )$ is without loops and $2$-cycles by the proof of \cite[thm.\ 0.6]{GHJ}.  Since $c_{ \cT }( u,\cU ) \neq 0$ we have $\Image \varphi \neq 0$.  Hence $\Image \varphi$ is indecomposable by Proposition \ref{psubobjects}, so we can set $m = \Image \varphi$.  
\end{proof}

\subsection{Dimension vectors of indecomposables are positive $c$-vectors of $\cC( Z )$}

In this subsection we prove an inclusion in the opposite direction to the one in the previous 
subsection, namely that the dimension vectors of all indecomposable objects of $\mod( \cT(T) )$ are positive 
$c$-vectors of $\cC(Z)$.

\bpr{dimension_vectors_of_indecomposables_are_c-vectors}
Let $\cT \subseteq \cC( Z )$ be a cluster tilting subcategory and let $m \in \mod( \cT )$ be indecomposable.

There exists a cluster tilting subcategory $\cU \subseteq \cC( Z )$ and an indecomposable $u \in \cU$ such that
\[
  c_{ \cT }( u,\cU ) = \dim_{ \cT }( m ). 
\]
\epr

\begin{proof}
There is a diagonal $v = \{ v_0,v_1 \}$ such that $\Hom_{ \cC( Z ) }( -,\Sigma v )\big|_{ \cT } = m$.  Set
\begin{equation}
\label{equ:forcingU}
  a_0 = v_1^+ \;,\; a_1 = v_0^+ \;,\;
  b_0 = v_0^- \;,\; b_1 = v_1^-
\end{equation}
and $u = \{ b_0,b_1 \}$, $u^{ \ast } = \{ a_0^-,a_1^- \}$, see Figure \ref{fig:dimension_vectors_of_indecomposables_are_c-vectors2}.
\begin{figure}
  \centering
    \begin{tikzpicture}[scale=2.5]
      \draw (0,0) circle (1cm);

      \draw[very thick] ([shift=(-80:1cm)]0,0) arc (-80:95:1cm);      
      \draw[very thick] ([shift=(125:1cm)]0,0) arc (125:250:1cm);            

      \draw (-80:0.97cm) -- (-80:1.03cm);
      \draw (-70:1.14cm) node{$a_0 = v_1^+$};
      \draw (95:0.97cm) -- (95:1.03cm);
      \draw (84:1.15cm) node{$b_0 = v_0^-$};
      \draw (110:0.97cm) -- (110:1.03cm);
      \draw (109:1.14cm) node{$v_0$};      
      \draw (125:0.97cm) -- (125:1.03cm);
      \draw (132:1.26cm) node{$a_1 = v_0^+$};
      \draw (250:0.97cm) -- (250:1.03cm);
      \draw (241:1.20cm) node{$b_1 = v_1^-$};
      \draw (265:0.97cm) -- (265:1.03cm);
      \draw (265:1.11cm) node{$v_1$};

      
%
                        
      \draw[blue] (110:1cm) .. controls (110:0.45cm) and (265:0.45cm) .. (265:1cm);

      \draw (95:1cm) .. controls (95:0.65cm) and (250:0.65cm) .. (250:1cm);

      \draw (0,0.6cm) node{$u$};

      \draw[blue] (0,-0.6cm) node{$v$};
      \draw (0.235,-0.569cm) node{$=u^{ \ast }$};      

    \end{tikzpicture} 
  \caption{Diagonals for the proof of Proposition \ref{pdimension_vectors_of_indecomposables_are_c-vectors}.}
\label{fig:dimension_vectors_of_indecomposables_are_c-vectors2}
\end{figure}
Note that in fact, $u = \Sigma v$, $u^{ \ast } = v$.    Let $U$ be a triangulation such that $u$ is in $U$ with mutation $u^{ \ast }$, and let $\cU$ be the cluster tilting subcategory corresponding to $U$.  

Let $T$ be the triangulation corresponding to $\cT$.    Since $u$ and $u^{ \ast }$ cross, there is a non-zero morphism $u \xrightarrow{ \delta } \Sigma u^{ \ast }$ and an induced morphism $\varphi$ as in Equation \eqref{equ:phi}.  The diagonal $v$ crosses at least one diagonal in $T$ since $m \neq 0$, so $\llbracket a_0,b_0 \rrbracket$ is connected to $\llbracket a_1,b_1 \rrbracket$ by at least one diagonal in $T$.  Proposition \ref{psubobjects}(2) says $\Image \varphi = \Hom_{ \cC( Z ) }( -,\Sigma h )\big|_{ \cT }$ where $h = \{ h_0,h_1 \}$ and the $h_i$ are provided by Lemma \ref{lU_improved}.  With the choices in Equation \eqref{equ:forcingU} we must have $h_0 = v_0$, $h_1 = v_1$, so $h = v$ and
\[
  \Image \varphi = \Hom_{ \cC( Z ) }( -,\Sigma v )\big|_{ \cT } = m.
\]

Since $m \neq 0$ we have $\Image \varphi \neq 0$, so $\delta \not\in [\Sigma \cT]$.  Hence Theorem \ref{t2}(2) gives
\[
  c_{ \cT }( u,\cU ) = \dim_{ \cT }\Image \varphi.
\]
The theorem applies because the quiver of each cluster tilting subcategory of $\cC( Z )$ is without loops and $2$-cycles by \cite[thm.\ 0.6]{GHJ}.  Combining the last two displayed equations proves the proposition.
\end{proof}


\subsection{Positive $c$-vectors of $\cC(Z)$ and dimension vectors of arcs of $Z$}

\bde{virtual_dimension_vectors}
If $e \neq f$ are in $\overline{ Z }$ then $\{ e,f \}$ is a {\em virtual arc}.  If $t = \{ t_0,t_1 \}$ is a diagonal, that is, if $t_0,t_1$ are non-neighbouring elements of $Z$, then $\{ e,f \}$ {\em crosses} $t$ if $e < t_0 < f < t_1 < e$ or $e < t_1 < f < t_0 < e$.

Let $T$ be a triangulation.  The {\em dimension vector of the virtual arc $\{ e,f \}$}, denoted
\[
  \dim_{ \cT( T ) }\big( \{ e,f \} \big) \in \Ksp\big( \cT( T ) \big)^{ \astsmall },
\]
is defined by
\begin{equation}
\label{d-vect}
  \langle \dim_{ \cT( T ) }\big( \{ e,f \} \big),t \rangle
  =
  \left\{
    \begin{array}{ll}
      1 & \mbox{if $\{ e,f \}$ crosses $t$,} \\
      0 & \mbox{otherwise}
    \end{array}
  \right.
\end{equation}
when $t \in T$ is a diagonal.

If $e \neq f$ are in $Z$, then the arc $\{ e,f \}$ is also a virtual arc so has a dimension vector.  The set of non-zero dimension vectors of arcs is denoted
\[
  D_{ \cT( T ) } =
  \big\{ \dim_{ \cT( T ) }\big( \{ e,f \} \big) \,\big|\, e \neq f \mbox{ are in } Z \big\} \setminus 0 .
\]
\ede

Our next theorem gives an explicit description of the set of positive $c$-vectors of $\cC(Z)$ in terms of the dimension vectors of the arcs of $Z$.
\bth{c-vectors_and_virtual_dimension_vectors}
If $T$ is a triangulation of $Z$, then
\[
  C^+_{ \cT( T ) }\big( \cC( Z ) \big) = D_{ \cT( T ) }.
\]
\eth

\begin{proof}
Observe that if $e \neq f$ are in $Z$ and we set $v = \{ e,f \}$, then 
\[
  \dim_{ \cT( T ) }\big( \{ e,f \} \big)
  = \dim_{ \cT( T ) }\big( \Hom_{ \cC( Z ) }( -,\Sigma v )\big|_{ \cT( T ) } \big).
\]
Hence $D_{ \cT( T ) }$ is the set of non-zero dimension vectors of indecomposable objects of $\mod( \cT )$, so the theorem follows from Propositions \ref{pc-vectors_are_dimension_vectors_of_indecomposables} and \ref{pdimension_vectors_of_indecomposables_are_c-vectors}.
\end{proof}

\sectionnew{Decompositions of the sets of positive $c$-vectors of $\cC(Z)$ via Borel subalgebras of $\sl_\infty$}
\lb{Borel}
Fix a triangulation $T$ of $Z$ and set
\[
C^+:= C^+_{\cT(T)}(\cC(Z)) \subset \Ksp(\cT(T))^*.
\]
In this section we classify the ideals of $C^+$ and prove that each of them can be identified with the roots of a Borel subalgebra 
of $\sl_\infty$. We also classify the Borel subalgebras of 
$\sl_\infty$ that appear in this way. These results are contained in Theorem \ref{tClassif-R+} and \ref{tXef-Borel}.
\subsection{The ideals of $C^+$}
\label{elemR}
\bde{support}
For a cluster tilting subcategory $\cT$ of a 2-CY category $\cC$, define the support of $c \in \Ksp( \cT )^{ \astsmall }$ by 
\[
  \supp\,c = \{ t \in \cT \,|\, \mbox{$t$ is indecomposable with $\langle c,[t] \rangle \neq 0$} \}.
\]
\ede
\bre{R_for_CZ}
\thref{c-vectors_and_virtual_dimension_vectors} implies that the ideals are $C^+$ are precisely the 
nonempty subsets $X \subseteq D_{ \cT( T )}$ which satisfy the following:
\begin{enumerate}
\setlength\itemsep{4pt}

  \item  If $d,d' \in D_{ \cT(T) }$ have $\supp\,d \supseteq \supp\,d'$, then $d \in X \Rightarrow d' \in X$.
  
  \item  If $d_1,\ldots,d_n \in X$ are given then there exists a $d \in X$ with $\supp\,d \supseteq \supp\,d_i$ for each $i$.

\end{enumerate}  
\ere
\bde{Xef}
For $e,f \in \ol{Z}$ we set
\[
  X_{ e,f } = \Big\{ c \in C^+ \Big| \supp\,c \subseteq \supp\,\dim_{ \cT( T ) }\big( \{ e,f \} \big) \Big\}.
\]
Note that
\[
  X_{ e,f } = \Big\{ d \in D_{ \cT( T ) } \Big| \supp\,d \subseteq \supp\,\dim_{ \cT( T ) }\big( \{ e,f \} \big) \Big\}
\]
by \thref{c-vectors_and_virtual_dimension_vectors}.
\ede
\bth{Classif-R+} For every cluster tilting subcategory $\cT(T)$ of the category $\cC(Z)$, the ideals of $C^+$ are precisely 
the nonempty subsets $X_{e, f} \subseteq C^+$ for $e, f \in \ol{Z}.$
\eth
We prove the theorem in two steps. In \S \ref{1Classif-R+} we show that every ideal of $C^+$ must be of the form 
$X_{e,f}$ for some $e, f \in \ol{Z}$. In \S \ref{Xef-str} we give an explicit description of the sets $X_{e, f }$ for all $e, f \in \ol{Z}$ 
and deduce from it that each nonempty $X_{e, f }$ is an ideal of $C^+$.

In \thref{Xef-Borel} we construct an explicit additive identification of each of the sets $X_{e,f}$ 
with the roots of a Borel subalgebra of $\sl_\infty$ or $\sl_n$. 
\subsection{Relation between ideals $C^+$ and the sets $X_{ e.f }$}
\label{1Classif-R+}
\bpr{each_X_is_Xef}
Each ideal $X$ of $C^+$ has the form $X = X_{ e,f }$ for certain $e \neq f$ in $\overline{ Z }$.
\epr

\begin{proof}
The proof will use the characterisation of ideals of $C^+$ from Remark \ref{rR_for_CZ}.

First, suppose that $X$ contains an element $d_{\rm max}$ which satisfies
\begin{equation}
\label{equ:c_max}
  d \in X \Rightarrow \supp\,d \subseteq \supp\,d_{\rm max}.	
\end{equation}
Then it is clear that
\[
  X = \Big\{ d \in D_{ \cT( T ) } \Big| \supp\,d \subseteq \supp\,d_{\rm max} \Big\}.
\]
But $d_{\rm max} = \dim_{ \cT( T ) }\big( \{ e,f \} \big)$ for certain $e,f \in Z$ and substituting this into the displayed formula shows $X = X_{ e,f }$.  

Secondly, suppose that $X$ does not contain an element $d_{\rm max}$ which satisfies Equation \eqref{equ:c_max}.  The set $Z$ is countable because it is a discrete subset of $S^1$.  Hence $D_{ \cT( T ) }$ is countable, so $X$ is countable.  Using Remark \ref{rR_for_CZ}, part (2), we can construct a sequence $\{ d_j \}_{ j \in \BN }$ in $X$ such that 
\begin{align}
\label{equ:seq1}
  & \supp\,d_1 \subseteq \supp\,d_2 \subseteq \supp\,d_3 \subseteq \cdots, \\
\label{equ:seq2}
  & \mbox{Each $d \in X$ has $\supp\,d \subseteq \supp\,d_j$ for some $j$,}
\end{align}
and it is clear that 
\begin{equation}
\label{equ:8.5}
  X = \Big\{ d \in D_{ \cT( T ) } \Big| \supp\,d \subseteq \supp\,d_j \mbox{ for some $j$} \Big\}.
\end{equation}
We have $d_j = \dim_{ \cT( T ) }\big( \{ e_j,f_j \} \big)$ for certain $e_j,f_j \in Z$.  By passing to a subsequence we can assume $e_j \rightarrow e$ and $f_j \rightarrow f$ for certain $e,f \in \overline{Z}$.  We will prove $X = X_{ e,f }$.  To prepare, suppose the diagonal $v = \{ v_0,v_1 \}$ has non-zero dimension vector $d = \dim_{ \cT( T ) }\big( \{ v_0,v_1 \} \big)$.  Then $v$ crosses at least one diagonal in $T$, so Lemma \ref{lU_improved_special_case} provides 
\[
  i_0 \leqslant s_0 < v_0 < i_1 \leqslant s_1 < v_1 < i_0
\]
in $Z$  such that $\{ i_0,v_1,s_1 \}$ and $\{ i_1,v_0,s_0 \}$ are triangles of $T$, see Figure \ref{fig:U_improved_special_case}.  Observe that if $g,h \in \overline{Z}$, then the implication
\begin{equation}
\label{equ:explanation_11a}
  \mbox{$t$ crosses $\{ v_0,v_1 \}$ $\Rightarrow$ $t$ crosses $\{ g,h \}$}
\end{equation}
for $t \in T$ is equivalent to 
\begin{equation}
\label{equ:explanation_11b}
  \mbox{$s_0 < g < i_1$ and $s_1 < h < i_0$, or vice versa.}
\end{equation}

The inclusion $X \subseteq X_{ e,f }$:  Suppose $d \in X$, that is $\supp\,d \subseteq \supp\,d_n$ for some $n$.  By Equation \eqref{equ:seq1} this implies that $\supp\,d \subseteq \supp\,d_j$ for each $j \geqslant n$, that is, $\supp\,\dim_{ \cT( T ) }\big( \{ v_0,v_1 \} \big) \subseteq \supp\,\dim_{ \cT( T ) }\big( \{ e_j,f_j \} \big)$ for each $j \geqslant n$.  Hence 
\[
  \mbox{$t$ crosses $\{ v_0,v_1 \}$ $\Rightarrow$ $t$ crosses $\{ e_j,f_j \}$}
\]
for $t \in T$ and $j \geqslant n$.  By \eqref{equ:explanation_11a} and \eqref{equ:explanation_11b} we can assume $s_0 < e_j < i_1$ and $s_1 < f_j < i_0$ for $j \geqslant n$.  This implies $s_0 < e < i_1$ and $s_1 < f < i_0$ whence $e \neq f$ and
\[
  \mbox{$t$ crosses $\{ v_0,v_1 \}$ $\Rightarrow$ $t$ crosses $\{ e,f \}$}
\]
for $t \in T$ by \eqref{equ:explanation_11a} and \eqref{equ:explanation_11b}.  Hence
\[
  \supp\,d = \supp\,\dim_{ \cT( T ) }\big( \{ v_0,v_1 \} \big) \subseteq \supp\,\dim_{ \cT( T ) }\big( \{ e,f \} \big)
\]
so $d \in X_{ e,f }$ as desired.  

The inclusion $X \supseteq X_{ e,f }$:  Suppose $d \in X_{ e,f }$, that is 
\[
\supp\,d = \supp\,\dim_{ \cT( T ) }\big( \{ v_0,v_1 \} \big) \subseteq \supp\,\dim_{ \cT( T ) }\big( \{ e,f \} \big).
\]
Hence 
\[
  \mbox{$t$ crosses $\{ v_0,v_1 \}$ $\Rightarrow$ $t$ crosses $\{ e,f \}$}
\]
for $t \in T$.  By \eqref{equ:explanation_11a} and \eqref{equ:explanation_11b} we can assume $s_0 < e < i_1$ and $s_1 < f < i_0$.  This implies that there is a $j$ for which $s_0 < e_j < i_1$ and $s_1 < f_j < i_0$, and by \eqref{equ:explanation_11a} and \eqref{equ:explanation_11b} this implies
\[
  \mbox{$t$ crosses $\{ v_0,v_1 \}$ $\Rightarrow$ $t$ crosses $\{ e_j,f_j \}$}
\]
for $t \in T$.  Hence $\supp\,\dim_{ \cT( T ) }\big( \{ v_0,v_1 \} \big) \subseteq \supp\,\dim_{ \cT( T ) }\big( \{ e_j,f_j \} \big)$, that is, $\supp\,d \subseteq \supp\,d_j$, so Equation \eqref{equ:8.5} implies $d \in X$ as desired.
\end{proof}
\subsection{The sets $X_{e,f}$ and $Y$}
\label{Xef-Y}
For the rest of this section we fix 
\[
e, f \in \ol{Z}
\]
such that $X_{e, f } \neq \varnothing$.
We will think of the diagonal $\{e,f \}$ as a straight line segment joining $e$ and $f$, and of crossing of diagonals as set theoretical intersection of such segments.
We will identify the open segment $\{e, f \}$ with the interval $(0,1)$, and define  
\begin{equation}
\label{Yset}
Y:= T \cap \{e, f\} \subset \{ e, f \} \cong (0,1).  
\end{equation}
\bpr{Y} (1) The set $Y \subset (0,1)$ has the following properties: 
\begin{enumerate}
\item[(a)] No point of $Y$ is a limit point of $Y$.  
\item[(b)] Every limit point of $Y$ in $(0,1)$ is both a left and a right limit point.
\end{enumerate}

(2) Every totally ordered set $Y$ which can be embedded into $(0,1)$ so that conditions (a)-(b) hold
arises from an Igusa--Todorov subset $Z \subset S^1$,
a triangulation $T$ of $Z$, and a pair $e, f \in \ol{Z}$. 
\epr
\begin{proof} (1) To show that $Y$ satisfies ($a$), assume that $y \in Y$ is the limit point of the sequence of distinct points $y_1, y_2, \ldots \in Y$. 
This means that there are diagonals
\[
\{g,h\}, \{g_1, h_1\}, \ldots \in T \quad \mbox{such that} \quad y = \{g, h\} \cap \{e, f \}, y_n = \{g, h \} \cap \{e_n,f_n \}. 
\]
By passing to a subsequence, we can assume that $\{g_n\}$ and $\{h_n\}$ converge. After a relabelling of $g$ and $h$ if needed, 
we get that
\[
\lim_{n \to \infty} g_n = g \quad \mbox{and} \quad \lim_{n \to \infty} h_n = h. 
\]
Since no point of $Z$ is a limit point of $Z$, $g_n=g$ and $h_n=h$ for sufficiently large $n$. Hence, $y_n =y$ for sufficiently large $n$, which is a contradiction. This proves that $Y$ has property ($a$). 

To verify that $Y$ has property ($b$), assume that $y$ is a limit point of $Y$ in the interval $(0,1)$ identified with $\{e,f\}$. 
This means that there exist
\[
\{g_1, h_1\}, \{ g_2, h_2 \} \ldots \in T \quad \mbox{such that} \quad y_n := \{g_n, h_n \} \cap \{e,f \} \to y \in \{ e, f \} 
\]
and $\{y_n\}$ are distinct.
By passing to a subsequence, we obtain that $\lim_{n \to \infty} g_n = g$ and $\lim_{n \to \infty} h_n = h$ for some 
$g \neq h \in \ol{Z}$, because $y \neq e$ and $y \neq f$.
Since $T$ is a triangulation, either $g_n =g$ for sufficiently large $n$ and the $h_n$ 
are distinct, or the other way around. Proceeding with the first case, we see that there should be a fountain in $T$ with limit point 
$h \in \ol{Z}$. The base of this fountain must be $g$. This means that there exits two sequences $\{h'_n\}$ and $\{h''_n\}$ in $Z$ such that
\[
h'_n \nearrow h, \; \; h''_n \searrow h \quad \mbox{and} \quad \{h'_n, g \}, \; \; \{ h''_n, g \} \in T.
\]
It follows that the sequences
\[
\{g, h'_n \} \cap \{ e, f \} \quad \mbox{and} \quad \{g, h''_n\} \cap \{e, f \} \; \; \mbox{in} \; \; Y
\]
left and right converge to $y$. 

(2) The statement of the second part is easy to check in the case $|Y| < \infty$; we leave the details to the reader. Let $Y \subset (0,1)$ be an infinite set with the properties $(a)$ and $(b)$. 

Choose $e \neq f \in S^1$ and denote by $\alpha$ and $\beta$ the two open semicircles of $S^1$ with end points $e$ and $f$. 
Identify $(0,1) \cong \{ e, f \}$. Choose a point $g \in \alpha$ and a subset $Z' \subset \beta$ such that
\[
Y = \{ y \in \{e, f \} \mid y = \{g,h\} \cap \{e, f\} \; \; \mbox{for some} \; \; h \in Z' \}.
\]
\begin{itemize}
\item[(i)] If $0 \cong e$ is a limit point of $Y$, choose a sequence of points $g'_n$ in a small neighborhood of $e$ in $\alpha$ converging to $e$. 
\item[(ii)] If $1 \cong f$ is a limit point of $Y$, choose a sequence of points $g''_n$ in a small neighborhood of $f$ in $\alpha$ converging to $f$. 
\end{itemize}
Let 
\[
Z := Z' \cup \{g\} \cup Z_0 \cup Z_1,
\]
where $Z_0 = \{ g'_n \mid n \in \Nset \}$ if $0 \in \ol{Y}$ and $Z_0 = \{ e\}$ otherwise; 
$Z_1 = \{ g''_n \mid n \in \Nset \}$ if $1 \in \ol{Y}$ and $Z_1 = \{ f\}$ otherwise. Let $T$ be the fan triangulation of $Z$ with base $g$. 
It is straightforward to check that $Z$ satisfies the Igusa--Todorov conditions, $T$ is a triangulation, and that 
$Y$ is the set associated to $X_{e,f}$.
\end{proof}
We will say that a set $Y$ is at most countable if it is finite or countable. The following proposition shows that 
the totally ordered sets that appear in \prref{Y} are precisely those that are at most countable and sequential
in the terminology from \S \ref{1.3}. 

\bpr{Y-descrp} For a totally ordered set $Y$, the following three conditions are equivalent:
\begin{enumerate}
\item $Y$ is at most countable and sequantial.
\item $Y$ is sequential and can be realized as a subset of $(0,1)$ in such a way that 
every limit point of $Y$ in $(0,1)$ is both left and right limit point.
\item $Y$ can be realized as a subset of $(0,1)$ in such a way that 
\begin{enumerate}
\item[(a)] no point of $Y$ is a limit point of $Y$ and
\item[(b)] every limit point of $Y$ in $(0,1)$ is both left and right limit point.
\end{enumerate}
\item $Y \cong Y_1 \sqcup Y_2 \sqcup Y_3$, where 
\[
\mbox{$Y_1 \cong \Nset$ or $\varnothing$},
\mbox{$Y_2 \cong X \times \Zset$ for an at most countable  totally ordered set $X$}, 
 \mbox{$Y_3 \cong - \Nset$ or $\varnothing$},
\]
with the lexicographic order on the Cartesian product and the relation $y_i < y_j$ for all $y_i \in Y_i$, $y_j \in Y_j$, $i <j$.
\end{enumerate}
\epr
\begin{proof} ($2 \Leftrightarrow 3$)
It is straightforward to show that the conditions in (2) and (3) are equivalent. 

($2 \Rightarrow 1$)
A subset $Y \subset (0,1)$, 
satisfying the conditions (2), is at most countable because for every $y \in Y$, the interval $(y,y^+)$ contains a rational number. 

($1 \Rightarrow 2$)
Assume that $Y$ is a totally ordered set satisfying (1). Using the assumption that it is at most countable, one constructs an embedding of $Y$ into
$(0,1)$. Without loss of generality we can assume that, if $Y$ has no least element, then $0$ is a limit point of $Y$, and 
if $Y$ has no greatest element, then $1$ is a limit point of $Y$.
Condition (a) in (3) is satisfied because $Y$ is sequential. To ensure that condition (b) is satisfied, one needs to modify the embedding.
Denote by $Y_{\lim}^l$ the set of left limit points of $Y$ without $\inf Y$ if it exists, and by $Y_{\lim}^r$ the set of right limit points of $Y$
without $\sup Y$ if it exists. We have a bijection
\[
b : Y_{\lim}^l \to Y_{\lim}^r \quad \mbox{given by} \quad b(y) := \min \{ y' \in Y_{\lim}^l \mid y' \geq y \}.
\]
Clearly, 
\[
Y \subseteq (0,1) \backslash \Big( \cup_{y \in Y_{\lim}^l} [y, b(y)] \Big).
\]
Choose a piecewise linear map $f : (0,1) \to (0,1)$ such that $f|_{[y,b(y)]} = \const$ for $y \in Y_{\lim}^l$ and 
$f$ is increasing in the rest of its domain. Composing the embedding $Y \subset (0,1)$ with the map $f$ gives an 
embedding that has the property (b) in (3).

($1 \Leftrightarrow 4$) Clearly, (4) implies (1). Now, let $Y$ be a totally ordered set which is at most countable and sequential. Consider the equivalence relation 
on $Y$ defined by $y \sim u$ if $y$ is obtained from $u$ by a finite chain of immediate predecessors or successors. There is a totally ordered set 
$X' := Y/\sim$. Let $X$ be the subset of $X'$ parametrizing the equivalence classes of $\sim$ that are isomorphic to $\Zset$. One easily sees that $X' = X_1 \sqcup X \sqcup X_3$ 
where $|X_1| \leq 1$, $|X_3| \leq 1$, and an easy reconstruction argument for $Y$ from $X'$ gives the desired form of $Y$ stated in (4). 

Note that the equivalence of (1) and (4) holds in greater generality when the condition of ``at most countable'' is removed from the statements of both parts.
\end{proof}
\bre{Zequiv} Analogously to the proof of \prref{Y-descrp}, one shows that for a cyclically ordered set $Z$ the following are equivalent:
\begin{enumerate}
\item $Z$ is at most countable, and every element of $Z$ has an immediate predecessor and immediate successor.
\item $Z$ can be realized as a subset of $S^1$ that satisfies the Igusa--Todorov conditions.
\item $Z \cong X \times \Zset$ for a cyclically ordered set $X$ which is at most countable with the lexicographic order on the Cartesian product.
\end{enumerate} 
\ere
\subsection{Structure of the sets $X_{e,f}$}
\label{Xef-str}
For a totally ordered set $Y$ and $y'\leq y'' \in Y$, set
\[
[y', y''] := \{ y \in Y \mid y' \leq y \leq y'' \} \subset Y.
\]
Denote the interval subset 
\[
\intt(Y) := \{ [y', y''] \mid y' \leq y'' \in Y \} \subset 2^Y.
\]
\bpr{Xef} Assume that $e, f \in \ol{Z}$ are such that $X_{e,f} \neq \varnothing$. Let $Y \subset \{e, f \} \cong (0,1)$ be given by \eqref{Yset}. The map
$\phi : X_{e,f} \to \intt(Y)$, given by
\[
\phi \Big( \dim_{ \cT( T ) }(\{a, b\}) \Big) := \big( \supp  \dim_{ \cT( T ) } \{a, b\} \big) \cap \{e, f \}
\]
for $a, b \in Z$ such that $\supp \dim_{ \cT( T ) } \{ a, b \} \subseteq \supp \dim_{ \cT( T ) } \{e, f \}$, 
is a bijection.
\epr
The proposition has the following corollary, which completes the proof of \thref{Classif-R+}.
\bco{Xef} For all $e, f \in \ol{Z}$ such that $X_{e, f} \neq \varnothing$, $X_{e, f}$ is an ideal of $C^+$.
\eco
\begin{proof}[Proof of \prref{Xef}] (1) First we show that $\phi( \dim_{ \cT( T ) }(\{a, b\})) \in \intt(Y)$ for all 
$a, b \in Z$ such that $\supp \{ a, b \} \subseteq \supp \{e, f \}$. Denote
\begin{align*}
&g' := \max \{ x \in \llbracket a^+, b^- \rrbracket \mid \{ a, x \} \in T\}, & g'' := \min \{ x \in \llbracket a^+, b^- \rrbracket \mid \{ b, x \} \in T\},
\\
& h' := \min \{ x \in \llbracket b^+, a^- \rrbracket \mid \{ a, x \} \in T\}, & h'' := \max \{ x \in \llbracket b^+, a^- \rrbracket \mid \{ b,  x \} \in T\},
\end{align*}
where the maximal and minimal elements exist due to \leref{U}. With respect to the cyclic 
order on $Z$ we have 
\[
a< g' < g'' < b < h'' < h' < a
\]
and the diagonals $\{ g', h' \}$ and $\{ g'', h'' \}$ belong to $T$. It follows that
\[
\dim_{ \cT( T ) } ( \{ a, b \} )
= \sum_{ g \in \llbracket g', g'' \rrbracket, \: h \in \llbracket h'', h' \rrbracket }  [ \{g, h \} ]^* 
\]
and from here that
\[
\phi \big( \dim_{ \cT( T ) } ( \{ a, b \} ) \Big) = [y', y'']
\]
where $y'$ and $y''$ are the smaller and bigger among the two points
\[
\{ g', h' \} \cap \{ e, f \}, \; \; \{ g'', h'' \} \cap \{ e, f \} \in Y \subset \{e, f \} \ \cong (0,1).
\]
The two intersections are nonempty because of the assumption $\supp \{ a, b \} \subseteq \supp \{e, f \}$.
This proves that the map $\phi : X_{e,f} \to \intt(Y)$ is well defined.

(2) Clearly the map $\phi$ is injective. So, it remains to prove that it is surjective.
Let $y' \leq y'' \in Y$. This means that there exist $g', g'', h', h'' \in Z$ such that 
\begin{equation}
\label{long<}
e < g' < g'' < f < h'' < h' < e
\end{equation}
in the cyclic order on $Z$ and $y'= \{ g', h' \} \cap \{ e, f \}$, $y''= \{ g'', h'' \} \cap \{ e, f \}$.
Denote by $\{ z',g',h' \}$ and $\{ z'',h'',g'' \}$ the triangles of $T$ such that 
$e$ and $z'$ are on the same side of $\{g', h' \}$, and $f$ and $z''$ are on the same side of $\{ g'', h'' \}$, 
guaranteed to exist by \leref{triangles}. We have
\[
\supp \dim_{ \cT( T ) } ( \{z',  z'' \} ) = \{ \{ g, h \} \mid g \in \llbracket g', g'' \rrbracket, h \in \llbracket h'', h' \rrbracket \} \in X_{e,f}, 
\]
because of \eqref{long<}. Hence, 
\[
\phi \big(  \dim_{ \cT( T ) } ( \{ z', z'' \} ) \big) =  \big\{ \{ g, h \} \cap \{ e, f \} \mid g \in \llbracket g', g'' \rrbracket, h \in \llbracket h'', h' \rrbracket \big\} =
[y', y''].
\]
This completes the proof of the proposition. 
\end{proof}
\subsection{The structure of the sets $X_{e,f}$ in terms of Borel subalgebras of $\sl_\infty$ and $\sl_n$}
\label{Xef-to-Y}
Recall the definition \eqref{extended} of the extentition $Y_{\ext}$ of a totally ordered set $Y$. 

\ble{Yiso} If $Y$ is a sequential, countable totally ordered set, then
\[
Y \cong Y_{\ext}
\] 
as totally ordered sets.
\ele
\begin{proof} We only need to consider the case when $Y$ has a least element. Denote it by $y_0$.
Denote by $y_1:=y_0^+, y_2:=y_1^+, \ldots$ the sequence of immediate successors starting from $y_0$. They exist because $Y$ is assumed to be infinite.
Thus, 
\[
Y = \{ y_0 < y_1 < \ldots \} \sqcup Y'
\]
with $y' > y_n$ for all $y' \in Y'$ and $n \in \Nset$. The map $Y \to Y_{\ext}$ given by 
\[
y_0 \mt -\infty, \; y_{n+1} \mt y_n, \; y' \mt y' 
\]
for $n \in \Nset$, $y' \in Y'$ is obviously an isomorphism of totally ordered sets.
\end{proof}
For $e, f \in \ol{Z}$ denote
\[
\Ksp(\cT(T))^*_{e, f} : = \Span X_{e, f} \subset \Ksp(\cT(T))^*.
\]
We will denote 
\[
\De^+ ( Y_{\ext} ) := 
\begin{cases}
\De^+_{\sl_\infty}(Y_{\ext}), & \mbox{if} \; \; |Y_{\ext}| = \infty,
\\
\De^+_{\sl_n}, & \mbox{if} \; \; |Y_{\ext}| = n
\end{cases}
\]
using the notation \eqref{De-Y}. Recall \eqref{Z-Y}. 

\bth{Xef-Borel} For every $e, f \in \ol{Z}$ such that $X_{e,f} \neq \varnothing$, there exists an isomorphism 
\[
\psi : \Ksp(\cT(T))^*_{e, f} \stackrel{\cong}{\lra} \Zset^{(Y_{\ext})}_0
\]
which restricts to a bijection
\[
X_{e, f} \stackrel{\cong}{\lra} \De^+(Y_{\ext}). 
\]
\eth
The right hand side of the bijection is the set of roots of a Borel subalgebra of $\sl_n$ or a  
Borel subalgebra of $\sl_\infty$, corresponding to a sequential (countable) totally ordered set $Y$.
\prref{Y}(2) implies that for each such Borel subalgebra, there exists an Igusa-Todorov subset $Z \subset S^1$, 
a triangulation $T$ of $Z$ and $e, f \in \ol{Z}$ such that $X_{e,f}$ is additively isomorphic to the 
set of roots of this Borel subalgebra.

\thref{Xef-Borel} follows from \prref{Xef} and the following lemma. For a totally ordered set $Y$, set
\[
\Zset^Y := \prod_{y \in Y} \Zset \ep_y \quad \mbox{and} \quad
1_{[y', y]} := \sum_{ y' \leq z \leq y} \ep_z \in \Zset^Y \quad \mbox{for} \quad y' \leq y \in Y.
\]
\ble{ZY-to-ZY0} Assume that $Y$ is an at most countable, sequential totally ordered set.
There is an isomorphism $\psi$ from 
\[
\Span \{ 1_{ [y, y']} \mid y \leq y' \in Y \} \subset \Zset^Y \quad \mbox{to} \quad
\Zset_0^{(Y_{\ext} )},
\]
satisfying 
\[
\psi(1_{[y', y]})= \ep_{y'} - \ep_{y^-}
\]
for all $y \leq y' \in Y$.
\ele
\begin{proof} Since $Y$ is sequential, every element of $\Span \{ 1_{ [y, y']} \mid y \leq y' \in Y \}$ can be uniquely represented in the 
form 
\begin{equation}
\label{els}
\sum_{i=1}^m c_i 1_{[y_i, z_i]}
\end{equation}
with $c_i \in \Zset \backslash 0$, $m\geq 0$, $y_i \leq z_i \in Y$, where
\begin{enumerate}
\item $z_i < y_{i+1}$ and 
\item $c_i \neq c_{i+1}$ if $z_i$ is the immediate predecessor of $y_{i+1}$. 
\end{enumerate}
Define a map  $\psi : \Span \{ 1_{ [y, y']} \mid y \leq y' \in Y \} \to \Zset_0^{(Y_{\ext} )}$ by
\[
\psi \big( \sum_{i=1}^m c_i 1_{[y_i, z_i]} \big) := \sum_{i=1}^m c_i (\ep_{z_i} - \ep_{y_i^-} )
\]
for the elements \eqref{els}. For $x,y,z \in Z$ with $x \leq y$ and $y^+ \leq z$, we have
\[
\psi( 1_{[x,y]} ) + \psi( 1_{[y^+, z]}) = (\ep_y - \ep_{x^-}) + (\ep_z - \ep_y) = \ep_z - \ep_{x^-} = 
\psi( 1_{[x,z]} ).
\]
It is easy to see that this implies that $\psi$ is a homomorphism. The facts that $\ker \psi =0$ and $\Im \psi = \Zset_0^{(Y_{\ext} )}$
are straightforward from the definition of $\psi$.
\end{proof}

\sectionnew{The maximal ideals of $C^+$}
\label{max-R+}
As in the previous section, fix a triangulation $T$ of $Z$ and set
\[
C^+:= C^+_{\cT(T)}(\cC(Z)) \subset \Ksp(\cT(T))^*.
\]
In this section we classify the maximal ideals of $C^+$ and prove that $C^+$ has a unique maximal ideal if and only if the dual 
quiver of $T$ has no cycles. This completes the proof of \thref{3}.
\subsection{A classification of the maximal ideals of $C^+$}
\label{max-elem-R}
\bpr{maxR} Let $e, f \in \ol{Z}$. The set $X_{e,f}$ is a maximal ideal of $C^+$ if and only if each of the 
points $e, f \in \ol{Z}$ is either an ear of $T$ or a limit point of a leapfrog of $T$.
\epr
\begin{proof} First we prove that, if each of the points $e, f \in \ol{Z}$ is either an ear of $T$ or a limit point of a leapfrog of $T$, 
then $X_{e,f}$ is a maximal ideal of $C^+$. The different cases for $e$ and $f$ are analogous. We will consider 
only the case when $e$ is the limit point of a leapfrog of $T$ and $f$ is an ear of $T$. 

Assume that $X$ is an ideal of $C^+$ such that $X_{e,f} \subseteq X$. It follows from \thref{Classif-R+} that 
$X = X_{a,b}$ for some $a, b \in T$. Since $\{f^-, f^+\} \in X_{e,f} \subset X_{a,b}$, the diagonal 
$\{e,f\}$ crosses the diagonal $\{a,b\}$. This implies that one of the end points of $\{a, b\}$ coincides with 
$f$. 

Furthermore, $X_{e,f}$ contains a sequence of diagonals of $T$ comprising the tail of the leapfrog with limit point $e$. 
Therefore, $X_{a,b}$ also contains those diagonals; that is $\{a,b\}$ crosses them. Hence, the other end point of $\{a,b\}$ 
coincides with $e$. So $\{a,b\} = \{e,f \}$, which means that $X_{e,f}$ is a maximal ideal of $C^+$.   

Now, let $X$ be a maximal ideal of $C^+$ and apply \thref{Classif-R+} to get that $X = X_{e,f}$ for
some $e,f \in \ol{Z}$. Assume that $e \in \ol{Z}$ is not an ear of $T$ or a limit point of a leapfrog of $T$.
Then, either (case 1) $e$ is the limit point of a fountain of $T$ or (case 2) $e \in Z$ and $e$ is not an ear of $T$.

(Case 1) Denote by $\{d, a_n\}$ and $\{d, b_n \}$ the diagonals of $T$ comprising the fountain with limit point $e$. 
After removing some of the diagonals of the fountain, we can assume that $d$ and all $a_n$ are on one side of the diagonal $\{e,f\}$ 
and all points $b_n$ are on the other. Then 
\[
\dim_{\cT(T)} (\{ e, f\}) \subset \dim_{\cT(T)} ( \{a_n, f\}) 
\]
and 
\[
\{ a_k, f\} \in \dim_{\cT(T)} (\{a_n, f\}) \backslash \dim_{\cT(T)} ( \{e, f\}) 
\]
for $k >n$. Therefore, $X_{e,f} \subsetneq X_{a_n, f}$ which contradicts the maximality of $X_{e,f}$. 

(Case 2) If $e \in Z$ and $e$ is not an ear of $T$, then there is a diagonal of $T$ with an end point $e$. This means that 
either $e < e^+ < h < f$ for some $h \in Z$ such that $\{e,h\} \in T$ or 
$f < h < e^- < e$ for some $h \in Z$ such that $\{e, h \} \in T$. In the two cases, respectively,
\[
X_{e, f} \subsetneq X_{e^+, f} \quad \mbox{and} \quad X_{e,f} \subsetneq X_{e^-,f},
\] 
which again contradicts the maximality of $X_{e,f}$. Hence, each of the points $e, f \in \ol{Z}$ must be 
either an ear of $T$ or a limit point of a leapfrog of $T$.
\end{proof}
\bre{root-syst} The combination of \prref{maxR}, \thref{Classif-R+}, and \S \ref{sec:decom}
implies that for every triangulation $T$ of $Z$,
\[
C^+_{\cT(T)}(\cC(Z)) = \bigcup X_{e,f}, 
\]
where the union ranges over $e, f \in \ol{Z}$ that are either an ear of $T$ or a limit point of a leapfrog of $T$.
By \thref{Xef-Borel}, each of the sets in the right hand side is additively isomorphic to 
\begin{enumerate}
\item the positive root system of one of the Lie algebras $\sl_n$ or
\item the root system of a Borel subalgebra of $\sl_\infty$ corresponding to a sequential (countable) totally ordered set $Y$. 
\end{enumerate}
It is easy to see that the argument of \prref{Y}(2) shows that for each of these positive roots systems $\De^+$, there exists 
an Igusa-Todorov set $Z \subset S^1$ and a triangulation $T$ of $Z$ such that 
\[
\mbox{$\De^+$ is additively isomorphic to $X_{e,f}$, as in \thref{Xef-Borel}},
\]
where $e$ and $f$ are ears or limit points of leapfrogs of $T$. The set $Z$ and the triangulation $T$ are constructed
explicitly from $Y$. 
\ere
Parts (1), (3) and (4) of \thref{3} follow from Theorems \ref{tClassif-R+} and \ref{tXef-Borel}, \prref{maxR} and \reref{root-syst}.
\subsection{Acyclic cluster tilting subcategories}
\label{acyclic}
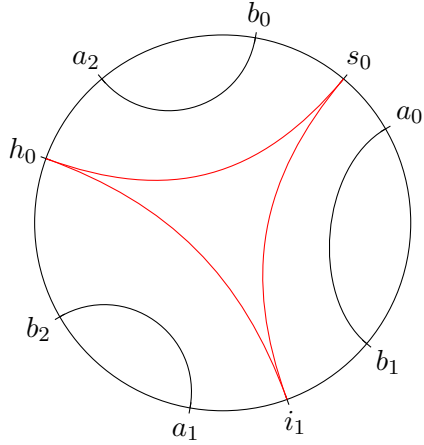
\begin{figure}
  \centering
    \begin{tikzpicture}[scale=2.5]
      \draw (0,0) circle (1cm);

      \draw (-40:0.97cm) -- (-40:1.03cm);
      \draw (-40:1.15cm) node{$b_1$};
      \draw (30:0.97cm) -- (30:1.03cm);
      \draw (30:1.15cm) node{$a_0$};
      \draw (50:0.97cm) -- (50:1.03cm);
      \draw (50:1.13cm) node{$s_0$};
      \draw (80:0.97cm) -- (80:1.03cm);
      \draw (80:1.13cm) node{$b_0$};
      \draw (130:0.97cm) -- (130:1.03cm);
      \draw (130:1.13cm) node{$a_2$};
      \draw (160:0.97cm) -- (160:1.03cm);
      \draw (160:1.13cm) node{$h_0$};
      \draw (210:0.97cm) -- (210:1.03cm);
      \draw (210:1.13cm) node{$b_2$};
      \draw (260:0.97cm) -- (260:1.03cm);
      \draw (260:1.13cm) node{$a_1$};
      \draw (290:0.97cm) -- (290:1.03cm);
      \draw (290:1.13cm) node{$i_1$};

      \draw (-40:1cm) .. controls (-40:0.6cm) and (30:0.6cm) .. (30:1cm);
      \draw (80:1cm) .. controls (80:0.6cm) and (130:0.6cm) .. (130:1cm);
      \draw (210:1cm) .. controls (210:0.6cm) and (260:0.6cm) .. (260:1cm);

      \draw[red] (50:1cm) .. controls (50:0.3cm) and (290:0.3cm) .. (290:1cm);
      \draw[red] (50:1cm) .. controls (50:0.3cm) and (160:0.3cm) .. (160:1cm);
      \draw[red] (290:1cm) .. controls (290:0.3cm) and (160:0.3cm) .. (160:1cm);



    \end{tikzpicture} 
  \caption{If a triangulation $T$ contains the black diagonals, then it also contains the red diagonals.  This forces non-acyclicity of $Q_T$.}
\label{fig:football}
\end{figure}

In this subsection we prove \thref{3}(2):
\bpr{QT-no-cycl} If the dual quiver $Q_T$ of the triangulation of $T$ has no cycles, then $C^+_{\cT(T)}(\cC(Z))$ has a unique maximal ideal
(equal to itself). 
\epr
\ble{football} Assume that $T$ is a triangulation of $Z$ which contains a configuration of black diagonals like the one on Figure \ref{fig:football}, where 
\[
a_2 \leq b_2 < a_1 \leq b_1 < a_0 \leq b_0 < a_2.
\]
Then $T$ contains a triangle consisting of diagonals; that is the dual graph $Q_T$ is not acyclic.
\ele
\begin{proof} Applying \leref{U_improved} to the vertices $a_1, b_1, a_0, b_0$ implies that there exists a triangle $\{i_1, h_0, s_0\}$ of $T$ such that 
\begin{equation}
\label{3-int}
i_1 \in \llbracket a_1, b_1 \rrbracket, \quad s_0 \in \llbracket a_0, b_0 \rrbracket, \quad h_0 \in \llbracket b_0, a_1 \rrbracket. 
\end{equation}
Since the arcs $\{i_1, h_0\}$ and $\{s_0, h_0\}$ do not intersect the diagonals  $\{ a_1, b_2\}$ and $\{ a_2, b_0 \}$, it follows that 
\[
h_0 \in \llbracket a_2, b_2 \rrbracket.
\]
This property and the first two properties in \eqref{3-int} imply that the arcs of the triangle $\{i_1, h_0, s_0\}$ of $T$ are diagonals.
\end{proof} 
\begin{proof}[Proof of \prref{QT-no-cycl}] Assume that $C^+_{\cT(T)}(\cC(Z))$ has more than one maximal ideal. It follows from \prref{maxR} that 
the total number of ears and limit points of leapfrogs of $T$  is at least 3. This implies that $T$ contains a configuration of diagonals like the one in Figure \ref{fig:football}. 
By \leref{football}, $Q_T$ is not acyclic, which is a contradiction.
\end{proof}
\medskip
\noindent
{\bf Acknowledgements.}
PJ was supported by EPSRC grant EP/P016014/1.  MY was supported by NSF grants DMS-1303038 and DMS-1601862, and by a Scheme 2 grant from the London Mathematical Society.

PJ thanks Louisiana State University and MY thanks Newcastle University for hospitality during visits.  We thank Ken Goodearl and Jon McCammond for helpful comments on ordered sets. We are grateful to the anonymous referee for numerous valuable suggestions.


\end{document}